\documentstyle{amsppt}

\NoBlackBoxes\nologo %
\hsize=15.5cm \vsize=20.0cm

\def\half{\frac{1}{2}}
\def\bt{\boxtimes}
\def\gr{{\text{gr}}}

\def\inv{^{-1}}
\def\?{{\bf{??}}}
\def\dgla{differential graded Lie algebra\ }
\def\dg{differential graded\ }

\def\C{\Bbb C}
\def\H{\Bbb H}
\def\P{\Bbb P}

\def\lb{\lambda}
\def\ord{\text{\rm ord}}
\def\sgn{\text{\rm sgn} }

\def\Der{\text{\rm Der} }
\def\Der{\text{\rm Der}}
\def\Spec{\text{\rm Spec} }
\def\ad{\text{\rm ad}}
\def\bp{\bar{\partial}}
\def\ls{\vskip.25in}

\def\Q{\Bbb Q}
\def\[{\big[}
\def\]{\big]}

\def\O{\Cal O}

\def\bt{\boxtimes}

\def\Sym{\text{\rm Sym}}
\def\id{\text{\rm id}}

\def\g{\frak g}

\def\gt{\tilde{\frak g}}

\def\gl{\frak{gl}}

\def\h{\frak h}

\def\m{\frak m}
\def\k{\frak k}
\def\endf{\frak { end}}
\def\hom{\frak {hom}}

\def\1/2{\frac{1}{2}}

\def\gsh{\frak g^\sharp}
\def\D{\frak D}
\def\U{\frak U}
\def\I{\Cal I}

\def\SS{\Cal S}

\def\im{\text{\rm im}}
\def\simto{\overset{\sim}\to{\to}}
\def\2{{[2]}}
\def\l{\ell}

\def\sbr #1.{^{[#1]}}
\def\scl #1.{^{\lceil#1\rceil}}
\def\spr #1.{^{(#1)}}
\def\emp #1. {{\it{#1\ }}}
\def\nl{\newline}
\def\bt{\boxtimes}
\def\inv{^{-1}}
\def\dgla{differential graded Lie algebra\ }
\def\dg{differential graded\ }

\def\C{\Bbb C}
\def\H{\Bbb H}
\def\P{\Bbb P}

\def\lb{\lambda}
\def\ord{\text{\rm ord}\ \ }
\def\sgn{\text{\rm sgn} }

\def\Der{\text{\rm Der} }
\def\Der{\text{\rm Der}}
\def\Spec{\text{\rm Spec} }
\def\ad{\text{\rm ad}}
\def\bp{\bar{\partial}}
\def\ls{\vskip.25in}

\def\Q{\Bbb Q}

\def\O{\Cal O}

\def\bt{\boxtimes}

\def\Sym{\text{Sym}}
\def\id{\text{id}}

\def\g{\frak g}

\def\gt{\tilde{\frak g}}

\def\gl{\frak{gl}}

\def\h{\frak h}

\def\m{\frak m}
\def\k{\frak k}
\def\endf{\frak { end}}
\def\hom{\frak {hom}}

\def\1/2{\frac{1}{2}}

\def\coker{\text {\rm coker}}
\def\Hom{\text {\rm Hom}}
\def\gsh{\frak g^\sharp}
\def\D{\frak D}
\def\U{\frak U}
\def\I{\Cal I}
\def\alt{^{\text{\rm alt}}}
\def\SS{\Cal S}

\def\simto{\overset{\sim}\to{\to}}

\def\Jsh{J^\sharp}
\def\at{{_\@}}
\def\bp{\partial}
\topmatter
\title {Lie Atoms and their deformations}\endtitle
\author{Z. Ran}\endauthor
\thanks{ Research supported in part by NSA grant
 H98230-05-1-0063; v.070104}\endthanks
\address University of California, Riverside\endaddress
\email ziv.ran\@ucr.edu\endemail\subjclass 14D15\endsubjclass
\abstract A Lie atom is essentially a pair of Lie algebras and its
deformation theory is that of a deformation with respect to the
first algebra, endowed with a trivialization with respect to the
second. Such deformations occur commonly in Algebraic Geometry, for
instance as deformations of subvarieties of a fixed ambient variety.
Here we study some basic notions related to Lie atoms, focussing
especially on their deformation theory, in particular the universal
deformation. We introduce Jacobi-Bernoulli cohomology, which yields
the deformation ring, and show that, under suitable hypotheses,
infinitesimal deformations are classified by certain Kodaira-Spencer
data.\endabstract
\endtopmatter
\document Many deformation-theoretic problems and results in algebraic and complex geometry can be profitably formulated in terms of Lie algebras, more specifically {\it{differential graded Lie algebras}} or dglas. These problems includes, notably, the Kodaira-Spencer theory of deformations of complex structures. Nevertheless, there are fundamental deformation problems in geometry for which no Lie theoretic formulation is known. These include, notably, the deformation theory of submanifolds in a fixed ambient manifold, i.e. the local theory of the Hilbert scheme in algebraic geometry or the Douady space in complex-analytic geometry. A principal purpose of this paper is to remedy this situation.\par
To this end, and for what we consider its own intrinsic interest, we introduce and begin to study a notion which we
call {\it{Lie atom}} and which generalizes that of the (shifted)
quotient of a Lie algebra by a subalgebra (more precisely, a pair of
Lie algebras up to bracket-preserving quasi-isomorphism)). Actually,
it turns out to be preferable to work with a somewhat more general
algebraic object, consisting of a pair of Lie algebras $\g, \h^+$, a
Lie homomorphism $\g\to\h^+$, and a $\g-$module $\h\subset\h^+.$ A
special case of this is a {\it Lie pair}, where $\h=\h^+$.
Geometrically, a Lie atom can be used to control situations where a
geometric object is deformed while some aspect of the geometry
'stays the same' (i.e. is deformed in a trivialized manner);
specifically, the algebra $\g$ controls the deformation while the
module $\h$ and the algebra $\h^+$ control the trivialization. A
typical example of this situation is that of a submanifold $Y$ in an
ambient manifold $X$, where $\g$ is the Lie algebra of relative
vector fields (infinitesimal motions of $X$ leaving $Y$ invariant),
and $\h=\h^+$ is the algebra and $\g$-module of all ambient vector
fields, so that the associated Lie atom is just the shifted normal
bundle $N_{Y/X}[-1]$, and the associated deformation theory is that
of the Hilbert scheme or Douady space of submanifolds of $X$.\par
Our point of view is that a Lie atom possesses some of the formal
properties of Lie algebras. In particular, we shall see that there
is a deformation theory for Lie atoms, which generalizes the case of
Lie algebras and which in addition allows us to treat some
classical, and disparate,
 deformation problems. These include, on the one hand,
the Hilbert scheme, and on the other hand heat-equation
deformations, introduced in the first-order case by Welters [We].
 Here we will present a systematic development of
some of the rudiments of the deformation theory of Lie atoms, which
are closely analogous to those of (differential graded) Lie
algebras. See \cite{Rrel2} for an application of Lie atoms to the
so-called Knizhnik-Zamolodchikov -Hitchin connection on the moduli
space of curves. \par An essential tool in the deformation theory
of Lie atoms is the {\it Jacobi-Bernoulli complex}, a comultiplicative complex whose zeroth cohomology, dualized,
yields the {\it deformation ring} of the atom. This is a local ring which, in good cases, e.g. when the global automorphisms are trivial or 'nearly' so,
is the base of the universal deformation. The Jacobi-Bernoulli complex is an
analogue of the familiar Jacobi complex of a Lie algebra, but its
differentials require twisting by Bernoulli numbers (hence the
name), and this makes it slightly less than obvious that the square of the
differential vanishes. The proof of this vanishing necessitates brief excursions into the realms of Lie
identities and Bernoulli number identities, which occupy \S0. In
\S1 we develop some basic algebra on Lie atoms and the
Jacobi-Bernoulli complex, and introduce some fundamental examples.
In \S2 we give some definitions and remarks on deformation theory
for Lie atoms, and introduce the Kodaira-Spencer formalism. Finally
in \S3 we construct universal deformations under suitable hypotheses
of finiteness and automorphism-paucity. See the introductions to individual sections for additional background, motivation and more detailed descriptions. \heading 0.
Preliminaries\endheading\subheading{ 0.1 Lie identities} The purpose
of this subsection is to write down some elementary, but possibly
non-standard, identities involving iterated brackets in a Lie
algebra. These identities will be technically useful in what follows. First some notation. Let
$a,b$ be elements in a Lie algebra with bracket $[,]$ and set
$$a_\@b=[a,b]=-\ad(b)(a)$$ and inductively, for any natural number
$m$, $$a_\@b^m=(a_\@b^{m-1})_\@b=(-\ad(b))^m(a).$$ For a function
$f(a_1,a_2)$ with values in an abelian group, we set
$$f(a_1,a_2)\alt=f(a_1,a_2)-f(a_2,a_1).$$ Thus, for example
$$[a_1,a_2]\alt =2[a_1,a_2];$$ from the Jacobi identity, it is easy
to check that $$[a_1,{a_2}_\@b]\alt =[a_1,a_2]_\@b.$$ We need a
generalization of the latter formula: \proclaim{Lemma 0.1} We have,
for all $m\geq 0$:
$$[a_1,{a_2}_\@b^m]\alt=$$$$
\sum\limits_{i=0}^{[m/2]}(-1)^i
\left(\binom{m-i-1}{i}+2\binom{m-i-1}{i-1}\right)
[{a_1}_\@b^i,{a_2}_\@b^i]_\@b^{m-2i}\tag 0.1$$ where we set
$\binom{j}{-1}=0, j\geq 0.$
\endproclaim\demo{proof} The Jacobi identity yields
$$[{a_1},{a_2}_\@b^m]\alt=[{a_1},{a_2}_\@b^{m-1}]\alt\ _\@b-
[{a_1}_\@b,{a_2}_\@b^{m-1}]\alt$$$$= [{a_1},{a_2}_\@b^{m-1}]\alt\
_\@b- [{a_1}_\@b,({a_2}_\@b)_\@b^{m-2}]\alt.$$ From this, it is easy
to see inductively that we may write
$$[a_1,{a_2}_\@b^m]\alt=
\sum\limits_{i=0}^{[m/2]}c_{i,m}
[{a_1}_\@b^i,{a_2}_\@b^i]_\@b^{m-2i}\tag 0.2$$ where the
coefficients $c_{i,m}$ satisfy  $$c_{0,m}=c_{0,m-1}, m\geq 2$$ so
that $c_{0,m}=1, m\geq 1, c_{0,0}=2$; and the recursion
$$c_{i,m}=c_{i,m-1}-c_{i-1,m-2}, 1\leq i\leq [m/2].\tag 0.3$$ To solve
this recursion set formally for $j\geq 1$ $$g_j(x)=\sum
c_{i,i+j}x^j.$$ Then we easily check that $c_{1,2}=-2$ to that
$g_1(x)=1-2x$, and (0.3)m translates into
$$g_j(x)=(1-x)g_{j-1}(x), j>1.$$ Thus $$g_j(x)=(1-2x)(1-x)^{j-1}$$
which yields (0.1).

\enddemo
\subheading{0.2 Bernoulli numbers} The Bernoulli numbers $B_n$ can
be defined by the generating function
$$C(x)=\sum\limits_{n=0}^\infty \frac{B_n}{n!}x^n=\frac{x}{e^x-1}
=-\frac{x}{2}+\frac{x}{2}\coth(\frac{x}{2}).$$ We set
$$c_n=\frac{B_n}{n!}, B(x)=C(x)+\frac{x}{2}=\frac{x}{2}\coth(\frac{x}{2}).$$
 Thus $c_0=1, c_1=-1/2, c_{2m+1}=0,\forall m>0.$ Moreover,
$$c_n=\sum\limits_{1\leq i\leq m\leq
n}(-1)^i\binom{m}{i}\frac{i^n}{(m+1)n!},\forall n\geq 1.\tag 0.4$$
$\lceil$ This formula will not be needed, but may be proved as
follows. Set $y=e^x-1$ so
$$C(x)=\frac{\log(1+y)}{y}=\sum\limits_{m=0}^\infty(-1)^my^m/(m+1).$$ The
binomial expansion yields
$$y^m=\sum\limits_{i=0}^m\sum\limits_{n=0
}^\infty(-1)^{i-m}\binom{m}{i}\frac{(ix)^n}{n!}.$$ Then the fact
that ord$_x(y^m)=m$, yields for $m\geq 1$
$$y^m=\sum\limits_{i=1}^m\sum\limits_{n\geq
m}(-1)^{i-m}\binom{m}{i}\frac{(ix)^n}{n!}.\rfloor$$ Now the freshman
calculus identity $d\coth x/dx=1-\coth^2x$ easily yields the
identities $$C^2=-xC'+(1-x)C, \tag 0.5$$ $$B^2=-xB'+B+x^2/4\tag
0.6$$ hence the quadratic recursion
$$(2m+1)c_{2m}=-\sum\limits_{i=1}^{m-1}c_{2i}c_{2m-2i}, m>1.$$
It is not hard to see from (0.6) that $B$ has the remarkable
property that the $\Q[x]$-module generated by all its derivatives is
closed under multiplication. We shall not need this fact as such,
but rather a precise form of a special case of it. First some
notation. Set $d=d/dx$ and
$$D_k=\frac{1}{k!}\prod\limits_{i=0}^{k-1}(-xd+k-i), \forall k\geq 1.\tag 0.7$$
Explicitly, in terms of power series,
$$D_k\sum a_ix^i=(-1)^k\sum\binom{i-1}{k}a_ix^i.$$
Multiplying $D_1B=-xB'+B$ by $B$ and using (0.6) to eliminate $B^2$ and its derivative $2BB'$, one can check easily that
$$B\cdot D_1B=D_2B.$$ We generalize this fact to higher derivatives
as follows \proclaim{Proposition 0.2} We have for $k\geq 1$
$$\matrix B\cdot D_kB=
\sum\limits_{i=0}^{[k/2]}c_{2i}x^{2i}D_{k+1-2i}B,\\
C\cdot D_kC=\sum\limits_{i=0}^{k}c_{i}x^{i}D_{k+1-i}C\endmatrix .\tag 0.8$$
Equivalently,  we have
$$\sum\binom{i-1}{k}c_ic_{m-i}=
-\sum\limits_{i=0}^{[k/2]} \binom{m-2i-1}{k+1-2i}c_{2i}c_{m-2i}
=-\sum\limits_{i=0}^{k} \binom{m-i-1}{k+1-i}c_{i}c_{m-i}
\tag
0.8bis$$
\endproclaim\demo{proof}
To start with, note that $D_kx=0, k\geq 1$, so $D_kB=D_kC$ and the two equations in (0.8) are equivalent; clearly the second equation is equivalent to (0.8bis).
Next,
note that $$D_\l B\equiv 1\mod x^\l,\forall
\l\geq 1,$$ hence $$x^iD_{k+1-i}B\equiv x^i\mod x^{k+1}.$$ Therefore to prove (0.8) it suffices to prove
\proclaim{$(*)_k$} $B\cdot D_kB$ is a constant linear combination of
$x^iD_{k+1-i}B, i\geq 0$
\endproclaim The coefficients of the linear combination are then
determined by examining the coefficients of $1, x, ...,x^k$. We will
prove $(*)_k$ by induction simultaneously with \proclaim{$(**)_k$}
$B^{k+1}$ is a constant linear combination of $x^{2i}D_{k-2i}B,
i\leq [\frac{k+1}{2}]$ where by definition $D_0B=B,
D_{-1}B=1$.\endproclaim \demo{proof} Firstly, assuming $(*)_k$ and $
(**)_k$, we deduce $(**)_{k+1}$ immediately by multiplying $(**)_k$
by $B$. So it remains to prove $(*)_{k+1}$. To this end we use $(**)_{k+1}$ to obtain an expression
$$B^{k+2}=\sum\limits_{i=0}^{[k/2+1]}a_ix^{2i}D_{k+1-2i}B\tag 0.9$$
for some constants $a_i.$ Comparing constant terms, it's clear that
$a_0=1.$ Now apply $(-xd+k+2)/(k+2)$  to (0.9). Using the operator
identity $$(-xd+r)x^m=x^m(-xd+r-m),$$ we get
$$\frac{1}{k+2}(-xd+k+2)B^{k+2}
=\sum\limits_{i=0}^{[k/2+1]}a'_ix^{2i}D_{k+2-2i}B\tag 0.10$$ with
$a'_0=1$. On the other hand, note that
$$\frac{1}{k+2}(-xd+k+2)B^{k+2}=B\frac{1}{k+1}(-xd+k+1)B^{k+1},$$
therefore, multiplying the analogue of (0.10) for $k+1$ by $B$
yields an expression for the same $\frac{1}{k+2}(-xd+k+2)B^{k+2}$ as
linear combination of the $Bx^{2i}D_{k+1-2i}B$. in which the $i=0$
term, i.e. $BD_{k+1}B$, appears with coefficient 1. Comparing the
two expressions and using $(*)_{k'}, k'\leq k$ now yields
$(*)_{k+1}.$\qed
\enddemo

\enddemo
It will be convenient to have a 'shifted' version of Proposition
0.2. For any integer $r$, define a shifted operator $D_k[r]$ by
$$D_k[r]=x^rD_kx^{-r}$$ or in terms of power series,
$$D_k[r]\sum a_ix^i=\sum \binom{i-1-r}{k}x^i.$$ The following
expansion can be verified easily
$$D_k[r]=\sum\limits_{j=0}^k\binom{r+1-j}{j}D_{k-j}.\tag 0.11$$
Combining this expansion with Proposition 0.2 and eq. (0.5), we conclude
\proclaim{Corollary 0.3} We have $$C\cdot D_k[r]C=
\sum\limits_{i=0}^{k}c_{i}x^{i}D_{k+1-i}[r]C-xC.\tag 0.12$$
Equivalently,
$$\sum\binom{i-1-r}{k}c_ic_{m-i}=
-\sum\limits_{i=0}^{k}
\binom{m-i-1-r}{k+1-i}c_{i}c_{m-i}-c_{m-1}\tag 0.12bis$$
\endproclaim Like (0.8bis), equation (0.12bis) is rather deceptive. Though its two sides look 'roughly' similar, they are in fact completely different in  nature, and that they happen to agree is a very special property of the Bernoulli function.
\comment
\par
Now set $$g=\coth(x).$$ We claim that for all $m\geq 0$, there exist
rational numbers $a_{i,m},b_{i,m}$ such that
$$g^m=\sum\limits_{i=-1}^{m-1}a_{i,m}g\spr i.,\tag ?$$
$$gg\spr m.=\sum\limits_{i=-1}^{m+1}b_{i,m}g\spr i.,\tag ?$$
where we set $g\spr -1. =1, a_{i,m}=0,\forall i\not\in [-1,m-1],
b_{i,m}=0,\forall i\not\in [-1,m+1].$ Moreover,
$$b_{m+1,m}=-1/(m+1), a_{m,m+1}=(-1)^m/m! \ ,$$ $$a_{-1,2m}=1,
a_{-1, 2m-1}=0, b_{-1,m}=0, m>0.$$ It is easy to see that this holds
for small $m$.. For the induction step, if $a_{*,m}, b_{*,m-1}$, it
is firstly clear that
$$a_{j,m+1}=\sum\limits_{i=j-1}^{m-1}b_{j,i}a_{i,m}\tag ??$$ Next,
differentiating ??  and solving for $g\spr m.$ yields $$g\spr
m.=(-1)^{m-1}(m-1)!(mg^{m-1}g'-\sum\limits_{i=1}^{m-1}a_{i-1,m}g\spr
i.)$$ Multiplying by $g$ and using the derivative of the $m+1$ case
of ?? yields $$gg\spr
m.=(-1)^{m-1}(m-1)!(\frac{m}{m+1}\sum\limits_{j=1}^{m+1}a_{j-1,m+1}g\spr
j.-\sum\limits_{i=1}^{m-1}\sum\limits_{j=-1}^{i+1}b_{j,i}a_{i-1,m}g\spr
j.)$$ That is, for $j\geq 0$,
$$b_{j,m}=(-1)^{m-1}(m-1)!(\frac{m}{m+1}\sum\limits_{i=0}^{m-1}b_{j-1,i}a_{i,m}
-\sum\limits_{i=1}^{m-1}\sum\limits_{i=1}^{m-1}b_{j,i}a_{i-1,m})\tag
??$$ It is easy to see inductively that $$i\equiv j\mod 2
\Rightarrow b_{i,j}=0.\tag ??$$ Now, inductively plugging the
equations ?? into themselves   we see that
$$\matrix b_{j,m}=&\\ (-1)^{m-1}&(m-1)!
(\frac{m}{m+1}\sum b_{j-1,i_1}b_{i_1,i_2}\cdots b_{i_m,-1}-\sum
b_{j,i_1}b_{i_1-1,i_2}b_{i_2,i_3}\cdots b_{i_m,-1})\endmatrix\tag ??
$$ where the indices satisfy $$i_t\leq m-t, \forall t$$ and $i_1\geq
0$ in the 1st sum, $i\geq 1$ in the second. Moreover in the first
sum, we automatically have $$i_{t-1}\leq i_t+1, t>1.$$ In the second
sum, the same inequality holds for $t>2$ while $i_1\leq i_2+2.$ Note
that in the first sum, if $i_t=m-t$ then $i_s=m-s,\forall s\geq t,$
hence
$$b_{i_s,i_{s+1}}=-1/(m-s), \forall s\geq t.$$ In the second sum,
the same thing holds if $t>1$. On the other hand if $t=1$, that is
$i_1=m-1$, then $m-3\leq i_2\leq m-2$ but since $b_{m-2,m-2}=0$ we
actually have $i_2=m-3$, hence
$$i_t=m-t-1,b_{i_t,i_{t+1}}=-1/(m-t-1), \forall t>1.$$ Now set
$$\matrix b_{j,m,\l,r}=&\\ (-1)^{m-1}&(m-1)!
(\frac{m}{m+1}\sum b_{j-1,i_1}b_{i_1,i_2}\cdots b_{i_\l,r}-\sum
b_{j,i_1}b_{i_1-1,i_2}b_{i_2,i_3}\cdots b_{i_\l,r})\endmatrix\tag ??
$$ both sums with similar constraints as above.; in particular,
$r\leq m-\l-1$. As noted above, for each monomial appearing in ??,
with indices $i_1,...,i_\l$, there is a uniquely determined number
$k\leq\l$ such that for all $t>k$ we have $i_t=m-t.$ Of course if
$r<m-\l-1$ then $k=\l$, therefore in this case
$$b_{j,m,\l,r}=b_{j,m-1,\l,r}.$$ Repeating the argument, we conclude
that $$b_{j,m,\l,r}=b_{j,\l+r+1,\l,r}.$$
\endcomment
\par \heading 1. Lie atoms: Basic notions\endheading
An inclusion $\g\to\h^+$ of Lie algebras, and more generally a homomorphism of dglas, constitutes in its own right a complex ('mapping cone') endowed with a bracket. Unfortunately, this complex is not usually a dgla (e.g. because the differential is not a derivation with respect to the bracket). Nevertheless, the structure involved is worth encoding, and this is accomplished through the notion of Lie atom.\par
 This section takes up the definition and initial study of Lie atoms. In \S1.1 we give the definition, some elementary remarks and constructions, and a basic list of standard examples, drawn mainly from geometry. Typically such examples involve a 'relative' situation, such as the inclusion of a submanifold in an ambient manifold. They will be used as a sort of 'benchmark' as we develop the theory.\par For technical reasons it is necessary to define a Lie atom in a slightly different, and finer, manner from the above naive notion, viz. as a $\g$-homomorphism $\g\to\h$ of a Lie algebra $\g$ to a $\g$-module. From such a homomorphism one can always construct a 'universal hull' $\h^+$, which is a Lie algebra receiving a Lie algebra homomorphism $\g\to\h^+$. This is why our definition is a refinement of the naive one.\par
 In \S1.2 we give the construction and basic properties of the Jacobi-Bernoulli complex associated to a Lie atom, which is a (nonobvious) extension of the Jacobi (standard) complex of a Lie algebra or dgla. This complex plays a fundamental role in the deformation theory of Lie atoms. The hardest part of the argument is the proof that the construction yields a complex, i.e. that the square of the differential vanishes. This proof requires the Bernoulli identities established in \S0.2. From the Jacobi-Bernoulli complex, we construct (see Theorem 1.2.1) the \emp{deformation ring}.  $R(\gsh)$ of a Lie atom $\gsh$, which will later be seen as the base of the universal deformation. We give a version (see Corollary 1.2.3) of the usual 'deformations minus obstructions' estimate on the dimension of the deformation ring, which is important in applications.  Then in \S1.3 and \S1.4 we touch on  the fundamental notions of {\it{atomic representation}} and {\it{universal enveloping atom}}, which are natural analogues of the corresponding notions for Lie algebras.
\subheading{1.1 Definition, examples, remarks} Unless otherwise
mentioned, Lie algebras will be understood over an arbitrary
commutative unitary ring $S$, which will usually be a
$\Q$-algebra.\proclaim{Definition 1.1.1} By a {\rm{ Lie atom}} (for
'algebra to module') we shall mean the data $\gsh$ consisting of
\item{(i)} a Lie algebra $\g$;
\item{(ii)} a $\g$-module
$\h$;
\item{(iii)} a $\g$-module homomorphism
$$i: \g\to\h,$$
where $\g$ is viewed as a $\g$-module via the adjoint
action.\nl\noindent
 If
$i$ is injective, $\gsh$ is said to be a {\text{ 'pure'}} Lie
atom.\nl If $\h$ is a Lie algebra and $i$ is a Lie homomorphism,
$\gsh$ is said to be a 'self-contained' Lie atom or a 'Lie pair'.
\endproclaim
\remark{Remarks 1.1.2}\par (1)
 Hypothesis (iii) means
explicitly that, writing $\langle  \ ,\ \rangle $ for the
$\g$-action on $\h$, we have
$$i([a,b])=\langle  a,i(b)\rangle =-\langle  b,i(a)\rangle .\tag 1.1.1$$
\comment(2) The definition and basic properties of Lie atoms (as in
this section) can be presented without the hypothesis $i$ injective;
let's call the resulting notion {\it{weak}} Lie atom. The
injectivity hypothesis will be important, however, when we discuss
the deformation theory of Lie atoms.\nl
\endcomment
\par (2) There is an obvious naive notion of atomic morphism of Lie
atoms, hence also of atomic isomorphism and quasi-isomorphism
(morphism inducing isomorphism on cohomology).  Of course one can
also talk about sheaves of Lie atoms, \dg\ Lie atoms, etc. Any Lie
atom is viewed as a complex in degrees 0,1, and we shall generally
consider two atoms to be equivalent if they are atomically
quasi-isomorphic. Accordingly, a {\it{morphism}} of Lie atoms would
be understood in the sense of the derived category, i.e. a homotopy
class of a composition of naive atomic morphisms and inverses of
naive atomic quasi-isomorphisms. Thus, for any Lie algebra $\g$, the
complex $\g\to 0$ yields a Lie atom equivalent to (and identified
with) $\g$; embedding $\g$ diagonally in $\g\oplus\g$ yields a Lie
atom equivalent to $\g[-1].$ For any Lie atom $(\g,\h)$, note the
'tautological' morphism $(\g,\h)\to \g.$\nl (3) Given a Lie atom
$\gsh$ as above, note that $\ker(i)$ is an ideal of $\g$ and $\gsh$
is an extension of the pure Lie atom $((\g/\ker(i))\to\h)$ by the
Lie algebra $\ker(i)$. Thus the notion of Lie atom is an
amalgamation of those of pure Lie atom and lie algebra.
\endremark A basic notion is that of a {\it hull} of a Lie atom.
Given a Lie atom $\gsh=(\g,\h,i),$ a { hull} for $\gsh$ is by
definition a Lie algebra $\h^+$ with a map $\h\to\h^+$ such that the
composite $\g\to\h^+$ is a Lie homomorphism and that the given
action of $\g$ on $\h$ extends via $i$ to a 'subalgabra' action of
$\g$ on $\h^+$, i.e. so that
$$\langle  a,v\rangle =[i(a), v],\ \forall a\in\g, v\in\h^+.$$
Note that any atom admits a universal hull $\h^\dag$, which is
simply the quotient of the free Lie algebra on $\h$ by the ideal
generated by elements of the form
$$[i(a), v]- \langle  a,v\rangle ,\ \ a\in\g, v\in\h$$
(note that the action of $\g$ on $\h$ extends to an action on
$\h^\dag$ by the 'derivation rule'). The basic identity (1.1.1)
shows that the map $\g\to\h^\dag$ induced by $i$ is a Lie
homomorphism.
\par In what follows,
we shall always understand a Lie atom $\gsh$ to come with a choice
of hull $\h^+$. If $\h$ itself is a hull, i.e. if $\gsh$ is a Lie
pair, we always choose $\h^+=\h$. On the contrary, if no hull is
specified, we take $\h^+=\h^\dag.$\par Now when $\g,\h^+$ are
nilpotent, $G=\exp(\g)\to H^+=\exp(\h^+)$ is a homomorphism of groups, but in general
$\exp(\h)\subset H^+$ is not $G$-invariant. For this reason we need
to consider what is essentially the tangent space to the $G$-orbit
of $\exp(\h)$. Denote by $\g\at\h^i$ the subgroup of $\h^+$
generated by all the $a\at b^i, a\in\g. b\in\h$, and set for
$m\in\Bbb N\cup\{\infty\},$
$$\h\sbr m.=\h+\sum\limits_{i< m+1}\g\at\h^i.\tag 1.1.2$$
 Also set
$$\g\sbr m.=(\g, \h\sbr m.)\subset\g^+=(\g,\h^+) .\tag 1.1.3$$
 Then it is easy to
see that $\h\sbr m.$ is a $\g$-module. More generally the adjoint
action yields a multi-pairing
$$\g\otimes\h\sbr m_1.\otimes\ldots\otimes\h\sbr m_k.\to\h\sbr
m_1+\ldots+m_k.\tag 1.1.4$$ Then from the Campbell-Hausdorff formula
we conclude \proclaim{Lemma 1.1.3} If $\g,\h^+$ are nilpotent, then
the subset $\exp(\h\sbr\infty.)\subset H^+$ is invariant under left
or right $G$-multiplication.\endproclaim \comment
 We
denote by $\h\sbr m.$ the subobject (e.g. subspace, subsheaf,
subcomplex...), automatically a $\g$-submodule, generated by
brackets of $\leq m$ elements of $\h$, so that we have an ascending
filtration (finite, if $\h^+$ is nilpotent), called the {\it
bracket} filtration:
$$\h\sbr 0.=i(\g)\subset\h=\h\sbr 1.\subset...\subset\h\sbr
m....\subset\h\sbr\infty.=\h^+.$$ This filtration is
bracket-compatible in the sense that
$$[\h\sbr i.,\h\sbr j.]\subset
\h\sbr i+j..$$\endcomment

\par \remark{ Stock Examples 1.1.4}
We present a list of basic examples to
be returned to repeatedly as we develop the theory of Lie atoms.
\par{\bf A.} The {\it{general linear atom}}. This essentially the universal example, of which every other is a special case. If $j:E_1\to E_2$ is
any linear map of vector spaces, let $\g=\g(j)$ be the
{\it{intertwining algebra}} of $j$, i.e. the Lie subalgebra
$$\g=\g(j)\subseteq\gl(E_1)\oplus\gl(E_2)$$ given by
$$\g =\{ (a_1,a_2)|j\circ a_1=a_2\circ j\} .$$
Thus $\g$ is the 'largest' algebra acting on $E_1$ and $E_2$ so that
$j$ is a $\g-$homomorphism. We define
$$\gl (E_1<E_2):=(\g, \gl (E_2), i_2),\tag 1.1.5$$
with $\ \ i_2(a_1,a_2)=a_2.\ \ $ (and, it goes without saying,
choice of hull as $\gl (E_2)$). Thus when $j$ is injective, so is
$i_2.$\par Next,  define
$$\gl (E_1>E_2):=(\g, \gl (E_1), i_1),\tag 1.1.6$$
with $\ \ i_1(a_1,a_2)=a_1. \ \ $ Thus when $j$ is surjective, $i_1$
is injective. These are Lie pairs.  The two notions are obviously
dual to each other, but since we do not assume $E_1, E_2$ are
finite-dimensional, dualising is not necessarily convenient.\par Finally, define $$\gl(E_1\vee E_2):=(\g, \gl(E_1)\oplus\gl(E_2), i_1\oplus i_2).\tag 1.1.7$$ \par In a
more global vein, we may consider a vector bundle homomorphism
$j:E_1\to E_2$ and define $\gl(E_1< E_2)$ and $\gl(E_1>E_2)$
similarly.
\par The foregoing construction
admits a useful generalization to the case of complexes (of vector
spaces, locally free sheaves or generally objects in an abelian
category). For any complex $(E^.,
\partial)$, we denote by $\gl(E^.)$ the 'internal hom' general linear
algebra of $E^.$, that is, the \dgla whose term in degree $i$ is
given by
$$\bigoplus_j\Hom(E^j, E^{j+i})$$ and whose differential is given
by (signed) commutator with $\partial$. When $E^.$ is an injective
resolution, i.e. a complex of injectives, acylic in a unique degree,
note that $\gl(E^.)$ is acyclic in negative degree, hence
quasi-isomorphic to a nonnegative complex. Given a morphism
$$j:E^._1\to E_2^.$$ of complexes, there is likewise an intertwining
\dgla $\gl(j)$; when $j$ is termwise injective (resp. surjective),
$\gl(j)$ is a subalgebra of $\gl(E_2^.)$ (resp. $\gl(E^._1)$). In
any event, there are $\gl(j)$-linear homomorphisms\nl
$$i_k:\gl(j)\to\gl(E^._k), k=1,2,$$ and we define differential
graded Lie pairs
$$\gl(E_1^.<E_2^.)=(\gl(j), \gl(E^._2),
i_2),$$$$\gl(E_1^.>E_2^.)=(\gl(j), \gl(E_1^.), i_1).$$ The first
(resp. second) definition is especially useful when $j$ is termwise
injective (resp. surjective). Thus, consider a short exact sequence,
say of coherent sheaves on a projective scheme $X$
$$0\to A\to B\to C\to 0.$$ As is well known, this sequence can be
resolved into a short exact sequence of complexes of locally free
coherent sheaves $$0\to E_1^.\to E_2^.\to E_3^.\to 0.$$ In fact we
may assume that each $E^i_j$ is a finite direct sum of line bundles
$\O_X(k)$ and that $E_2^i=E_1^i\oplus E_3^i, \forall i.$ If $X$ is
smooth (e.g. $X=\P^n$), the complexes $E_i^.$ may be assumed
bounded. It is easy to see that the Lie atoms $\gl(E_1^.<E_2^.),
\gl(E_2^.>E_3^.)$ are, up to quasi-isomorphism, independent of the
resolution, so we may set
$$\gl(A<B)=\gl(E_1^.<E_2^.), \gl(B>C)=\gl(E_2^.>E_3^.).$$ As we
shall see, these Lie atoms control the formal germ at $B\to C$ of
the Quot scheme of $B$ (a special case of which is the Hilbert
scheme).\par
 {\bf B.} If $i:\g_1\to\g_2$
is a
  homomorphism of Lie algebras, then
$$\gsh :=(\g_1,\g_2,i)$$
is a Lie pair. More generally, if $\h$ is any $\g_1$ submodule of
$\g_2$ containing $i(\g_1)$, then
$$\gsh :=(\g_1,\h,i)$$
is a Lie atom, whose hull will be taken as the Lie subalgebra of
$\g_2$ generated by $\h$.\par Note that a general Lie atom
$(\g,\h,i)$ is essentially of this type, modulo replacing $\g$ and
$\h$ by their images in the hull $\h^+$ (though the map $\h\to\h^+$
is not necessarily injective).

 {\bf C.} Let $E$ be an invertible
  locally free sheaf on a ringed space $X$
(such as a real or complex manifold), and let $\D^i(E)$ be the sheaf
of $i-$th order differential endomorphisms of $E, i\geq 0,$ and
$$\D^{\infty}(E)=\bigcup\limits_{i=0}^{\infty}\D^i(E).$$
Then $\g=\D^1(E)$ and $\D^\infty(E)$ are Lie algebra sheaves and
$\h=\D^2(E)$ is a $\g-$submodule of $\D^\infty(E)$ , giving rise to
a Lie atom $\gsh$ with hull $\D^\infty(E)$,  which will be called
the {\it{Heat atom}} of $E$ and denote by $\D^{1/2}(E)$. Note that
if $X$ is a manifold then $\gsh$ is quasi-isomorphic as a complex to
$\Sym^2(T_X)[-1]$ and $\h\sbr m.=\D^{m+1}(E)$.
\par {\bf D.} Let
$$Y\subset X$$ be an embedding of manifolds (real or complex). Let
$T_{X/Y}$ be the sheaf of vector fields on $X$ tangent to $Y$ along
$Y$. Then $T_{X/Y}$ is a sheaf of Lie algebras contained in its
module $T_X$, giving rise to a  Lie pair
$$N_{Y/X}[-1]= (T_{X/Y}\subset T_X),$$
which we call the {\it{normal atom}} to $Y$ in $X$. Notice that
$T_{X/Y}\to T_X$ is locally an isomorphism off $Y$, so replacing
$T_{X/Y}$ and $T_X$ by their sheaf-theoretic restrictions on $Y$
yields a Lie atom that is quasi-isomorphic to, and identifiable with
$N_{Y/X}[-1]$.\par More generally, let $f:Y\to X$ be an arbitrary
morphism (e.g. a holomorphic map of complex manifolds). Then we have
a Lie algebra sheaf on $Y$: $$T_{X/Y}=\ker (f\inv T_X\oplus T_Y\to
f^*T_X)$$$$ =\{(u,v)\in\Der(f\inv\O_X)\oplus \Der(\O_Y):f^*\circ
u=v\circ f^*\}.$$ This sheaf, sometimes called the sheaf of
'$f$-related vector fields' admits as modules both $f\inv T_X$ and
$T_Y$, giving rise to Lie pairs on $Y$:
$$N_{f/X}[-1]=(T_{X/Y}\to f\inv T_X),$$$$N_{f/Y}[-1]=(T_{X/Y}\to
T_Y).$$ $N_{f/X}[-1]$ is pure whenever $f$ is generically immersive.
Note that when $f$ is an inclusion of a submanifold, there is an
obvious quasi-isomorphism of Lie atoms $N_{Y/X}[-1]\to T_{f/X}[-1].$
\comment
 Via the graph of $f$, we get an embedding $Y\subset Y\times
X$, whence a Lie algebra $$T_f= T_{Y\times X/Y}$$ and a Lie atom
$$N_{f,Y,X}[-1]=(T_f, T_{Y\times X}). $$  If $f$ is generically
immersive, so that $df$ is generically injective, we get another Lie
atom
$$N_{f/X}[-1]=(T_f, p_X^*T_X).$$
\endcomment
 If $f$ is submersive,
the Lie atom $N_{f/X}[-1]$ is equivalent to the algebra of vertical
vector fields $T_{f}=\ker (T_Y\overset df\to\to f^*T_X)$ while
$N_{f/Y}[-1]$ is pure.

 {\bf E} In the situation of the
previous example with $Y\subset X$, let $\I_Y$ denote the ideal
sheaf of $Y$. Then $\I_Y.T_X$ is also a Lie subalgebra of $T_X$
giving rise to a  Lie pair
$$T_X\otimes\O_Y[-1]:=(\I_Y.T_X\subset T_X)\sim_{\text{qis}}
(T_{X/Y}\subset T_X\otimes T_Y).$$ Note that via the
embedding of $Y$ in $Y\times X$ as the graph of the inclusion
$Y\subset X$, $T_X\otimes\O_Y[-1]$ is quasi isomorphic as Lie atom
to $N_{Y/Y\times X}[-1]$, so this example is essentially a special
case of Example D.\qed
\endremark\ls \subheading{1.2 Jacobi-Bernoulli complex}\par Our purpose
here is to define the {\it{Jacobi-Bernoulli}} (or JacoBer, for
short) complex $\Jsh(\gsh)$ associated to a Lie atom
$\gsh=(\g\to\h)$, which is to play an analogous role in the
deformation theory of $\gsh$ as the Jacobi complex $J(\g)$ in the
deformation theory of a Lie algebra $\g$. We begin with the case
where $\gsh$ is a {\it Lie pair.}\par For a Lie pair $\gsh$,
$\Jsh(\gsh)$ is a complex in nonpositive degrees whose terms $K^.$
are defined as follows:$$K^0=\bigoplus\limits_{j>0}K^{0,j},$$
$$K^i=\bigoplus\limits_{j\geq 0}K^{i,j},$$where
$$K^{i,j}=\bigwedge\limits^{-i}\g\otimes\Sym^j\h.$$ Next we define the
differential $d^i:K^i\to K^{i+1}, i<0.$ This will be a direct sum
$$d^i=\bigoplus d^{i,j,j'}:K^{i,j}\to K^{i+1, j'}, j'\leq j+1.$$ In
particular, $K^.$ will {\it not} be a double complex, though it has
a natural increasing degree filtration $F_.K^.$ defined by
$$(F_rK)^i=\bigoplus\limits_{j\leq i+r}K^{i,j}.$$ It turns out that
the role of $F_.$ will be analogous to that of the 'stupid'
filtration on the ordinary Jacobi complex. To define the $d^i$, we
start with $d^{-1}.$ As in \S0.2,  we denote by $c_t=B_t/t!$
the normalized Bernoulli coefficient.Define
$$d^{-1, m, m-t+1}:\g\otimes\Sym^m\h\to\Sym^{m-t+1}\h,\ 0\leq t\leq m$$
by
$$a\otimes b^m\mapsto c_t(i(a)_\@b^t).b^{m-t}=c_t\sum\limits_{j=0}^{m-t}b^j(i(a)_\@b^t)b^{m-t-j}\tag 1.2.1$$  Also define, as
in the usual Jacobi complex,
$$d^{-2,-0,0}:\bigwedge\limits^2\g\to\g,$$
$$a_1\wedge a_2\mapsto [a_1,a_2].$$ This then determines the other
differentials via the 'derivation rule':
$$d^{-n, m,
m-t+1}:\bigwedge\limits^n\g\otimes\Sym^m\h\to
\bigwedge\limits^{n-1}\g\otimes\Sym^{m-t+1}\h,$$
$$a_1\wedge...\wedge a_n\otimes b^m\mapsto$$$$
c_t\sum\limits_{i=1}^n (-1)^{i-1}a_1\wedge...\widehat{a_i}...\wedge
a_n\otimes(i(a_i)_\@b^t).b^{m-t}{\hskip 2cm}\tag 1.2.2
$$$$+\sum\limits_{i<j}(-1)^{i-j-1}a_1\wedge...\wedge [a_i,a_j]\wedge...
\hat{a_j}...\wedge a_n\otimes b^m$$ 
It is not obvious that these differentials
define a complex, because of the twisting by Bernoulli numbers. We
summarize the essential properties of $\Jsh(\gsh)$ as follows
\proclaim{Theorem 1.2.1}(i)$(\Jsh, F.)$ is a functor from the
category of Lie atoms over $S$ to that of comultiplicative,
cocommmutative and coassociative filtered complexes over $S$.\par
(ii)The filtration $F_.$ is compatible with the comultiplication and
has associated graded
$$F_i/F_{i-1}=\bigwedge\limits^i(\gsh).$$\par
(iii) The quasi-isomorphism class of $\Jsh(\gsh)$ depends only on
the quasi-isomorphism class of $\gsh$ as Lie
atom.\endproclaim\demo{proof} The hard part is proving that $K^.$ as
above is a complex, i.e. that $d^{-n+1}d^{-n}=0.$ Using the
definition of $d^.$, one reduces easily, first to that case $n=2$,
then to the vanishing of the component of $d^{-1}d^{-2}$ going from
$K^{-2,m-1}=\bigwedge\limits^2\g\otimes\Sym^{m-1}\h$ to
$K^{0,1}=\h$, that is, to proving that $$\sum\limits_{i=0}^md^{-2,
m-1, i}d^{-1, i, 1}(a_1\wedge a_2\otimes b^{m-1})=0, \forall a_1, a_2\in
\g, b\in\h.\tag 1.2.3$$  Plugging into the definitions, (1.2.3) means
$$\sum\limits_{i+j\leq m-1}c_ic_{m-i}[a_1{ _\@}b^j, a_2 {_\@}b^i]\alt
\at b^{m-1-i-j}
+c_{m-1}[a_1,a_2]{_\@}b^{m-1}=0.\tag 1.2.4$$ 
We break the big summation in two subsums I and
II depending on whether $j\leq i$ or $j>i.$ Using Lemma 0.1, I can
be evaluated as $$\sum\limits_{j+2\alpha\leq i}[a_1\at
b^{j+\alpha},a_2\at b^{j+\alpha}]\at
b^{m-1-2j-2\alpha}\cdot$$$$\cdot(-1)^\alpha
(\binom{i-j-1-\alpha}{\alpha}+2\binom{i-j-1-\alpha}{\alpha-1})
c_ic_{m-i}.$$ (Note if $m$ is even, only even $i$'s appear.) Setting
$r=j+\alpha$, and using the elementary formula
$$\sum\limits_{\alpha=0}^x(-1)^\alpha\binom{y}{\alpha}=(-1)^x\binom{y-1}{x}$$
the coefficient of $[a_1\at b^r, a_2\at b^r]b^{m-1-2r}$ can be
evaluated as $$I_{r,m}=\sum\limits_i
(-1)^r(\binom{i-r-2}{r}-2\binom{i-r-2}{r-1})c_ic_{m-i}.$$ Referring
to \S0.2, the latter is none other than the degree-m term in
$$(D_r[r+1]B-2D_{r-1}[r+1]B)B.$$ The subsum II can be analyzed
in the same way. Now setting $r=i+\alpha$, we get $$II=\sum [a_1\at
b^r, a_2\at b^r]\at b^{m-1-2r}II_{r,m}$$ where
$$II_{r,m}=\sum\limits_i
(-1)^{r-i}(\binom{m-i-r-1}{r-i+1}-2\binom{m-i-r-1}{r-i})c_ic_{m-i}$$
As in \S0.2, this is just the degree$-m$ term in
 $$-\sum\limits_{i\leq r\ {\text{}}} c_ix^iD_{r+1-i}[r+1]C+
 +2\sum\limits_{i\leq r-1\ {\text{}}} c_ix^iD_{r-i}[r+1]C.$$
 By Corollary 0.3,
we conclude $I+II+c_{m-1}[a_1,a_2]_{\@}b^{m-1}=0$, 
completing the proof that $K$ is a complex. [A
more conceptual proof of this is much to be desired !!]\par Now
given that $\Jsh(\gsh)$ is a complex, functoriality and filtration
are obvious, as is the assertion about the graded. As for the
comultiplication, its definition is directly analogous to that of
the comultiplication in the ordinary Jacobi complex, as developed in
\cite{R}, based on the natural map
$$\bigwedge\limits^i\g\otimes\Sym^m\h\to
\bigoplus\limits_{\overset{i_1+i_2=i}\to{m_1+m_2=m}}
(\bigwedge\limits^{i_1}\g\otimes\Sym^{m_1}\h) \otimes
(\bigwedge\limits^{i_2}\g\otimes\Sym^{m_2}\h)$$
The proof of the required properties
of this comultiplication follows along the same lines as the proof
in \cite{R} of the analogous assertions for the Jacobi complex.
\par Finally, as for (iii), a quasi-isomorphism
$\gsh_1\to\gsh_2$ of Lie atoms induces quasi isomoprhisms
$$\bigwedge\limits^i(\gsh_1)\to\bigwedge\limits^i(\gsh_2)$$ and the
spectral sequence of a filtered complex then shows that the induced
map $\Jsh(\gsh_1)\to\Jsh(\gsh_2)$ is a quasi-isomorphism.\qed
\proclaim{Corollary 1.2.2} Assume $\gsh$ is acyclic in nonpositive degrees. Then there is a second-quadrant spectral sequence with $E_1$ term
$$E_1^{p,q}=\Sym^{-q-2p}\H^1(\gsh)\otimes\Sym^{q+p}\H^2(\gsh)$$ whose abutment has degree 0 part equal to $\bigoplus E^{p,-p}_\infty$ where $$ E^{p,-p}_\infty=\bigoplus {\text{gr}}^{-p}_{F_.}\H^0(\Jsh(\gsh)).$$\endproclaim
\demo{proof}
For any complex $(K^.,d)$, its {\it{truncation in degree}} $r$ is defined by $$\matrix (K^{r]})^i=& K^i, i<r,\\&\ker d^r, i=r,\\&0, i>r.
\endmatrix$$ We have $$\matrix H^i(K^{r]})=& H^i(K), i\leq r,\\&0, i>r.\endmatrix$$ Then , it is easy to see that $$\H^0(\Jsh((\gsh)^{2]}))=\H^0(\Jsh(\gsh)).$$ Indeed this simply follows from the fact that, for all $r$, $$d_r(E_r^{p,-p})\subset E_r^{p+r, -p-r-1}$$ where, in the case of a JacoBer complex of a Lie atom, the latter group
invloves only $\H^1$ and $\H^2$ of the atom.\par Therefore, replacing $\gsh$ by its truncation in degree 2, we may assume $\gsh$ is acyclic in degrees $>2$. Then the Corollary becomes simply the spectral sequence associated to the $F.$ filtration, as in Thorem 1.2.1(ii).
\qed\enddemo

\enddemo In view of Theorem 1.2.1, there is a natural $S$-algebra
structure on
$$R_S(\gsh):=S\oplus\Hom_S(\H^0(\Jsh(\gsh)), S)$$ and we call the
latter the {\it deformation algebra} or {\it ring} of the Lie pair
$\gsh$ over $S$. We now take $S=\C$, and denote $R_S(\gsh)$ simply by $R(\gsh)$. It is a local ring with maximal ideal $\m=\H^0(\Jsh(\gsh))^*$. Note that the $\m$-adic filtration on $R=R(\gsh)$ is dual to the
$F.$ (multiplicative-order) filtration on $\Jsh(\gsh)$. As $R$ will usually be 'tested' by mapping it to a finite-dimensional $\C$-algebra, our mail interest will be in the formal completion
$$\hat{R}=\varprojlim\limits_e R/\m^e.$$
Now, the spectral sequence of Corollary 1.2.1 yields some information on the associated graded of the $\m$-adic filtration on $R$, i.e. the tangent cone  $\gr^.R$. Set
$$V=\H^1(\gsh)^*, W=\H^2(\gsh)^*.$$
Thus, in the above spectral sequence,
$$E_1^{p,-p}=\Sym^{-p}V^*, E_1^{p, -p+1}=\Sym^{-p-1}V^*\otimes W^*.$$ By construction, the cone $\gr^.(R)$ is a quotient of the symmetric algebra $\Sym^.V$, and from the spectral sequence,
we will obtain a description of the corresponding homogeneous ideal $I^.=\bigoplus I^r<\Sym^.V$.
\par To this end,
note that the spectral sequence endows $W$ with a decreasing filtration
$$W=W^2\supseteq W^3=\im (d_1^{-2, 2})^\perp\supseteq...\supseteq W^{r+1}=\im(d_{r-1}^{-r, r})^\perp... $$ (recall that the image of $d_{r-1}^{-r, r}$ is defined modulo the image of $d_{r-2}^{-r+1, r-1}$). The dual of $d_r^{-r-1, r+1}$ yields a map $$c_r:W^r/W^{r+1}\to \Sym^rV.$$ Now, note that the differential $$d_1^{-r, r}:\Sym^rV^*\to \Sym^{r-2}V^*\otimes W^*$$ is compatible with the comultiplication $\Sym^rV^*\to \Sym^{r-2}V^*\otimes\Sym^2V^*.$ Therefore, $d_1^{-r, r}$ is just the map extended in the obvious way from
$d_1^{-2, 2}:\Sym^2V^*\to W^*.$ A similar property holds for the other differentials $d_s^{-r, r}$. Dualizing this fact, we conclude that $${\text {gr}}^rR=(E_\infty ^{-r, r})^*=
\coker\left(\bigoplus\limits_{i=0}^{r-2} W^{r-i}/W^{r-i+1}\otimes \Sym^iV\overset{(\tilde{c}_r,...,\tilde{c}_2)}\to{\longrightarrow} S^rV\right )$$ where gr refers to the $\m$-adic filtration and the map $\tilde{c}_{r-i}$ on the $i$th summand is the one extended from $c_{r-i}$. In other words, the kernel $I^r$ of the surjection $$\Sym^rV\to \m^r/\m^{r+1}$$ is the portion in degree $r$ of of homogeneous ideal in the symmetric algebra $\Sym^.V$ generated by the images of the $c_{r-i}, i=0,...,r-2.$ Thus $I$ has a set of generators corresponding to a $\C$-basis of $W$. In particular, if $V$ and $W$ are finite-dimensional, then $I$ is generated by $\dim W$ many elements and $\Sym^.V$ is a regular noetherian ring of Krull dimension $\dim V$; by Krull's Theorem, it follows that the Krull dimension of $\gr^.R$, measured e.g. as $1+$ the degree of the appropriate Hilbert polynomial
$$p(r)=\dim_\C(\gr^r(R)), r>>0,$$ is at least $\dim V-\dim W$. Since the completion $\hat{R}=\hat{R}(\gsh)=\varprojlim\limits_e R/\m^{e}$ is
a complete noetherian local ring of the same Krull dimension as the tangent cone, we have proven
\proclaim{Corollary 1.2.3} If $\gsh$ is acyclic in nonpositive degrees and $h^1(\gsh)$ and $h^2(\gsh)$ are finite, then $$\dim\hat{R}(\gsh)\geq h^1(\gsh)-h^2(\gsh).\qed$$\endproclaim
This result is a fundamental 'a priori estimate' on dimension. Given the relation between $R(\gsh)$ and deformations, to be established in \S3 below, it may be viewed as an \emp{existence}. statement, asserting, whenever $h^1(\gsh)-h^2(\gsh)>0$, the existence of a nontrivial deformation. Corollary 1.2.3 unifies and extends a number of 'folklore' results, some of which have had
important applications. As just one example, one could mention Mori's theory of rational curves (see \cite{Kol}).
\par
 Next, we will define the JacoBer complex $\Jsh(\gsh)$
for a general Lie atom $\gsh=(\g,\h)$ with hull $\g^+=(\g, \h^+).$
$\Jsh(\gsh)$ will be a certain subcomplex of $\Jsh(\g^+)$, endowed
with a sub-stupid filtration that depends on the submodule
$\h\subset\h^+$.\par Recall the 'adjoint' filtration $(\h\sbr\cdot
.)$ on $\h^+$ (cf. (1.1.2)). It naturally induces a filtration on
each $\Sym^i\h^+$: namely
$$\Sym^i(\h^+)\sbr j.\subset\Sym^i(\h^+)$$ is the subgroup generated
by all $b_1\cdots b_i$ such that  for some $(j_1,...,j_i)$, we have
$$b_k\in\h\sbr j_k.,\forall k,
{\text {and}}\  \sum j_k=j.$$ There is a similar induced filtration
on $\bigwedge\limits^r\g\otimes\Sym^i\h^+$, declaring that elements
of $\g$ have filtration level 0. It is evident that this extended
bracket filtration is compatible with the differentials on
$\Jsh(\g^+)$, thus yielding a filtration of $\Jsh(\g^+)$ by
subcomplexes $\Jsh(\g^+)\sbr j.$. We set
$$\Jsh_m(\gsh):=F_m(\Jsh(\gsh)):=
F_m(\Jsh(\g^+))\cap \Jsh(\g^+)\sbr m..$$ Then $\Jsh_m(\gsh)$
inherits a comultiplicative structure from $\Jsh(\gsh)$, and we set
$$R_m(\gsh)=S\oplus\Hom_S(\H^0(\Jsh_m(\gsh)),S),$$ which is then an
$S$-algebra called the $m$th deformation algebra of the Lie atom
$\gsh.$

\comment

 An important, though elementary, remark about Lie atoms is
that a Lie atom indeed constitutes a 'Lie object' in the category of
complexes, i.e. that the natural map $\bigwedge\limits^2\gsh\to\gsh$
extends to a complex, called a Jacobi complex
$$J(\gsh):...\bigwedge\limits^{i+1}\gsh\to
\bigwedge\limits^i\gsh...\to\gsh.$$ Here $\bigwedge\limits^i\gsh$ of
course means exterior power as complex. Thus $J(\gsh)$ is actually a
double
complex $$\matrix\cdots&S^3\h&&&&\\ &\uparrow&&&&\\ \cdots&\g\otimes S^2\h&\to& S^2\h&&\\
&\uparrow&&\uparrow&&\\
\cdots&\bigwedge\limits^2\g\otimes\h&\to&\g\otimes\h&\to&\h\\
&\uparrow&&\uparrow&&\uparrow\\
\cdots&\bigwedge\limits^3\g&\to&\bigwedge\limits^2\g&\to&\g
\endmatrix\tag 1.2$$
with up arrows defined by $$g_1\wedge\cdots\wedge g_r\otimes
h_1\cdots
h_s\mapsto\frac{1}{r}\sum\limits_{j=1}^r(-1)^jg_1\wedge\cdots\widehat{g_j}\cdots\wedge
g_r\otimes i(g_j)h_1\cdots h_s\tag 1.3$$ and right arrows defined by
$$g_1\wedge\cdots\wedge g_r\otimes h_1\cdots
h_s\mapsto$$$$\frac{1}{r!}\sum\limits_{\sigma\in\SS_r}\sum\limits_{j=1}^s
\sgn(\sigma)[g_{\sigma(1)},g_{\sigma(2)}]\wedge g_{\sigma(3)}\cdots
g_{\sigma(r)}h_1\cdots h_s- $$$$g_{\sigma(2)}\wedge\cdots \wedge
g_{\sigma(r)}h_1\cdots\langle  g_{\sigma(1)},h_j\rangle \cdots
h_s$$$$ +g_{\sigma(1)}\wedge g_{\sigma(3)}\cdots \wedge
g_{\sigma(r)}h_1\cdots\langle  g_{\sigma(2)},h_j\rangle \cdots
h_s\tag 1.4$$
 The $\g-$linearity of $i$, plus the
$\frac{1}{r}$ factor make each square commute (if $i$ were a
$g$-derivation rather than $\g$-linear, each square would commute
without the $\frac{1}{r}$ factor, which would yield a \dgla, a
decidedly distinct notion from Lie atom). We shall return to the
Jacobi complex in \S3, in connection with universal deformations.
For now we content ourselves with a formal remark.
\endcomment
By construction then, we have a chain of complexes with inclusion
maps
$$\Jsh_1(\gsh)=\g\sbr 1.\to...\to \Jsh_m(\gsh)\overset{i_m}\to{\to}
 \Jsh_{m+1}(\gsh)
\to...\to \Jsh(\gsh):=\Jsh_\infty(\gsh)\to\Jsh(\g^+)\tag 1.2.5$$
inducing a chain of rings and homomorphisms $$ R(\g^+)\to
R(\gsh)=R_\infty(\gsh)\to\cdots\to
R_{m+1}(\gsh)\overset{\eta_m}\to\to R_{m}(\gsh)\to\cdots .\tag
1.2.6$$ \proclaim{Lemma 1.2.4} Suppose $\g$ is a \dgla with $H^{\leq
0}(\g)=0$.\par (i) Let $\h$ be a dg $\g$ module in nonnegative degrees.
Then the map $\eta_m$ in (1.2.6) is surjective and its kernel
contains $\m_{R_{m+1}(\gsh)}^{m+1}$.\par (ii) For any injection
$\h_1\to\h_2$ of nonnegative $\g$-modules, the induced map
$$H^0(\Jsh_m(\g,\h_1))\to H^0(\Jsh_m(\g,\h_2))$$ is
injective.\endproclaim \demo{proof} It suffices to prove (ii). This
is a standard spectral sequence argument. Note that by choosing a
complement to $\partial\g^0$ in $\g^1$, we may replace $\g^.$ by a
sub-dgla $\k$ in strictly positive degrees that is quasi-isomorphic
to $\g$ and yields a quasi-isomorphic atom $(\k, \h)$. Then
replacing $\gsh$ by $(\k,\h)$, the mapping cone of $i_m$ can be
represented by a complex in nonnegative degrees, hence has no
negative cohomology, which implies our assertion.
\qed\enddemo \proclaim{Definition 1.2.5}A dg Lie atom $(\g,\h)$ is
said to be {\it positive} if it is isomorphic as Lie atom to one
$(\g',\h')$ where $\g'$ (resp. $\h'$) exists only in positive (resp.
nonnegative) degrees.\endproclaim Clearly, the conclusion of Lemma
1.2.2(i) applies to any positive dg Lie atom.\remark{Remark 1.2.6}
When $\gsh$ is a sheaf of (possibly dg) Lie atoms on a topological
space $X$, there are sheaf-theoretic analogues of the JacoBer
complexes $\Jsh(\gsh), \Jsh_m(\gsh)$, analogous to the
sheaf-theoretic Jacobi complex (cf. \cite{Rcid}). These are
complexes defined on certain subset spaces $X\langle m\rangle,
X\langle\infty\rangle$, where the ordinary tensor, exterior or
symmetric products $\otimes, \bigwedge\limits^i, \Sym^i$ are
replaced by their 'external' analogues $\bt, \lambda^i, \sigma^i.$
This construction works well when $X$ is Hausdorff, but not
otherwise (e.g. when $X$ is a scheme in the Zariski topology). When
$X$ is non-Hausdorff, a version of the subset spaces for
Grothendieck topologies still works. However, for any topological
space, one can always replace a sheaf of Lie atoms or dg Lie atoms
$\gsh=(\g,\h)$ by a suitable acyclic (soft or injective) resolution
$(\g^.,\h^.)$, and then $(\Gamma(\g^.),\Gamma(\h^.))$ is a dg Lie
atom whose cohomology is the same as the sheaf cohomology of $\gsh$.
In the sequel, we understand constructions like the JacoBer complex,
applied to a Lie atom sheaf $(\g,\h)$ to mean the non-sheafy version
applied to $(\Gamma(\g^.),\Gamma(\h^.))$.\endremark
 \subheading{1.3 Representations} Now given a Lie
atom $\gsh=(\g,\h,i)$, by a {\it{left $\gsh-$module}} or {\it{left
$\gsh-$representation}} we shall mean the data of a pair $(E_1,E_2)$
of $\g-$modules with an injective $\g-$ homomorphism $j:E_1\to E_2$,
together with an 'action rule'
$$\langle  \ \ \rangle :\h\times E_2\to E_2,$$
satisfying the compatibility condition (in which we have written
$\langle  \ \rangle $ for all the various action rules):
$$\langle  \langle  a,v\rangle ,x\rangle =
\langle  a,\langle  v,x\rangle \rangle - \langle  v,\langle
a,x\rangle \rangle ,\tag 1.3.1$$
$$\forall a\in\g , v\in\h ,x\in E_2,$$ such that the $\h$-action
extends to a Lie action of the hull $\h^+$ (if no hull is specified,
so $\h^+$ is the universal hull, the latter condition is redundant).
In other words, a left $\gsh-$ module is just a homomorphism of Lie
atoms

$$\gsh\to\gl (E_1<E_2).$$
Similarly, a right $\gsh-$module is defined as a homomorphism of Lie
atoms\nl $\gsh\to\gl (E_1>E_2).$

\remark{Examples 1.1.4, cont.}  Refer to the previous examples.\item
{\bf A.} These are the tautological examples: $\gl(E_1<E_2)$ and
$\gl(E_1>E_2)$ with $(E_1,E_2)$ as left (resp. right) module in the
two cases $j$ injective (resp. surjective).\item {\bf B.}For a Lie
atom $\gsh=(\g_1,\g_2,i)$, $\gsh$ itself is a left  $\gsh$-module,
called the {\it{adjoint representation}} while $(\gsh)^*=(\g_2^*,
\g_1^*, i^*), *=$ dual vector space, is a right $\gsh-$module called
the {\it{coadjoint representation}}.\item {\bf C.} In this case
$(E,E)$ is a left and right $\gsh-$module, called a {\it{Heat
module}}.\item {\bf D.} For an inclusion $Y\subset X$, the basic
left module is $(\I_Y,\O_X)$ which is quasi-isomorphic to $O_Y.$ Of
course we may replace $\I_Y$ and $\O_X$ by their topological
restrictions on $Y$. The basic right module of interest is
$(\O_X|_Y,\O_Y)$. For a general map $f:Y\to X$, the basic left
$N_{f/X}[-1]$-module of interest is $(f\inv\O_X\to \O_Y)$. For a
submersion $f: Y\to X$, the basic left $T_{f/Y}[-1]$-module is
$(f\inv\O_X\to \O_Y)$.\item {\bf E.} Realizing $T_X\otimes\O_Y[-1]$ as $T_{X/Y}\to T_X\oplus T_Y$, the natural right module of interest is $$O_X\to O_Y.$$
\endremark
\comment \remark{Remark} 'Theoretically', only the action of $\h$
going from $E_1$ to $E_2$ 'should' be necessary for a module.
However the action on all of $E_2$ is needed in proofs and satisfied
in the examples we have in mind, so we included it. The fact that we
require an extension to $E_2$ rather than the dual notion of a
'lifting' to $E_1$ has to do with the fact that in our examples of
interest the maps $i$ and $j$ are injective.
\endremark
\endcomment
\subheading{ 1.4 Universal enveloping atom}

\comment Given a Lie atom $\gsh=(\g,\h,i),$ a {\it hull} for $\gsh$
is by definition a Lie algebra $\h^+$ with a $\g-$map $\h\to\h^+$
such that the given action of $\g$ on $\h$ extends via $i$ to a
'subalgabra' action of $\g$ on $\h^+$, i.e. so that
$$\langle  a,v\rangle =[i(a), v],\ \forall a\in\g, v\in\h^+.\tag 1.6$$
Note that there is a universal hull $\h^\dag$, which is simply the
quotient of the free Lie algebra on $\h$ by the ideal generated by
elements of the form
$$[i(a), v]- \langle  a,v\rangle ,\ \ a\in\g, v\in\h$$
(note that the action of $\g$ on $\h$ extends to an action on
$\h^\dag$ by the 'derivation rule'). In view of the basic identity
(1.1) it follows that the map $\g\to\h^\dag$ induced by $i$ is a Lie
homomorphism. This shows in particular that, modulo replacing $\g$
and $\h$ by their images in $\h^\dag$, any Lie atom $\gsh$ is
essentially of the type of Example B above (though the map
$\h\to\h^\dag$ is not necessarily injective).\par In what follows,
we shall always understand a Lie atom to come with a choice of hull
$\h^+$. If $\h$ itself is a hull, we always choose $\h^+=\h$. If no
hull is explicitly specified, we take $\h^+=\h^\dag.$\par
\endcomment
 We observe
next that there is an natural notion of 'universal enveloping atom'
associated to a Lie atom $\gsh=(\g,\h,i)$. Indeed, denoting by
$\U(\g)=\C1\oplus\U^+(\g)$ the usual enveloping algebra, let
$\U(\g,\h)$ be the quotient of the $\U(\g)$-bimodule
$\U(\g)\otimes\h\otimes\U(\g)$ by the sub-bimodule generated by
elements of the form
$$a\otimes v\otimes 1-1\otimes v\otimes a-1\otimes\langle  a,v\rangle\otimes 1,
1\otimes i(a)\otimes b-a\otimes i(b)\otimes 1,$$
$$\forall a,b\in\g, v\in\h.$$
\noindent{\bf{Sorites }}\par {1.}\  $\U(\g,\h)$ is a
$\U(\g)-$bimodule .\par 2.\ The map $i$ extends to a bimodule
homomorphism
$$i^\sharp:\U(\g)\to\U(\g,\h).$$\par
{3.}\   $\U(\g,\h)$ is universal with respect to these
properties.\par 4.\ $\U(\g,\h)$ is generated by $\h$ as either right
or left
 $\U(\g)-$module. Moreover the image of $\U(\g,\h)$
 in $\U(\h^+)$ is precisely
 the (left, right or bi-)$\U(\g)$-submodule of $\U(\h^+)$
 generated by $\h$. \par 5. An atomic analogue of the
 Poincar\'e-Birkhoff-Witt  formula holds: it states that, as vector spaces,
 $$\U(\g, \h)\simeq \Sym^{>0}({\text {im}} (i))
 \oplus\Sym^{\geq 0}(\g)\otimes{\coker}(i)$$ with the two summands
 corresponding, respectively, to ${\text{im}}(i^\sharp)$ and
 $\coker(i^\sharp)$. Moreover, under the usual PBW isomorphism
 $$\U(\g)\simeq\Sym(\g),$$ the kernel of $i^\sharp$ corresponds to
 $\Sym(\g)\ker(i),$ i.e. the ideal in $\Sym(g)$ generated by
 $\ker(i)$.
\par Thus $$\U(\gsh):=(\U(\g),\U(\g,\h),i)$$ forms an 'associative
atom' which we call the {\it{universal enveloping atom}} associated
to $\gsh$. From (5), it is easy to see that a quasi-isomorphism
$\phi:\gsh_1\to\gsh_2$ yields a quasi-isomorphism
$\U(\phi):\U(\gsh_1)\to\U(\gsh_2)$.
  \remark{Examples 1.1.4, cont. }\item {\bf A.} It is
elementary that the universal enveloping algebra of the interwining
Lie algebra $\g$ is simply the interwining associative algebra
$$\U(\g)=\{(a_1,a_2)|j\circ a_1=a_2\circ j\},$$
and so the universal enveloping atom of $\gl (E_1<E_2)$ (resp. $\gl
(E_1>E_2)$ is just $(\U(\g), \endf (E_2), i)\ $ (resp. $ (\U(\g),
\endf (E_1), i) .$\item {\bf B. } In this case it is clear that
$\U(\g_1,\h)$ is just the sub $\U(\g_1)-$bimodule generated by $\h$.
\endremark
 \heading 2. Atomic deformation theory basics\endheading
 The purpose of this section is to develop the deformation theory associated to a sheaf of Lie atoms $\gsh=(\g\to\h^+)$ on a topological space. This is a generalization to a 'relative' setting of the familiar Kodaira-Spencer deformation theory for a sheaf $\g$ of Lie algebras, using the same point of view as that as developed in \cite{Rcid}. This theory, which will be reviewed extensively below, is essentially {\it{geometric}} in character. It starts from the notion of a $\g$-deformation,
 over (or parametrized by) a finite-dimensional local $\C$-algebra $S$ with maximal ideal $\m_S$, as a certain type of torsor $G^\phi$ with respect to the sheaf of nilpotent groups $G_S=\exp(\g\otimes\m_S)$. Such a torsor can be 'realized' on any $\g$-module $E$, yielding a deformation $E^\phi$ of $E$. \par
 To a $\g$-deformation one can associate, following the ideas of Kodaira-Spencer, a formal algebraic object known as Kodaira-Spencer class which, in 'raw' form is represented by an element $$\phi\in\g^1\otimes\m_S$$ satisfying the integrability or Maurer-Cartan equation $$\partial\phi=-\half [\phi,\phi].$$ A more abstract formulation, in terms of the cohomologically-defined 'deformation ring' $R(\g)$, can also be given (see below). In this formulation, the element $\phi$ is replaced by a (multiplicative) homomorphism $\epsilon=\epsilon(\phi):R(\g)\to S$, called the Kodaira-Spencer homomorphism of the deformation.\par To clarify our perspective, we note that some treatments of deformation theory  {\it{define}} a deformation formally via the notions of Kodaira-Spencer class and integrability equation. By contrast, our approach {\it{deduces}} these from the geometric definition. While the geometric and formal notions of deformation are closely related in both directions, the exact formulation of this relation, and especially its descent to suitable equivalence or cohomology classes, is somewhat subtle, and depends on appropriate hypotheses of finiteness and automorphism-freeness. Defining a deformation formally via Kodaira-Spencer class sweeps these issues under the rug.\par
 Now extending this deformation theory to the case of a Lie atom $\gsh=(\g\to\h^+)$ means, essentially, considering a $\g$-deformation, in the form of a torsor $G^\phi$, together with a trivialization of the induced $H^+$-torsor $(H^+)^\phi$. For simplicity we ignore in this introduction the role of the $\g$-submodule $\h\subset\h^+$.
 To a $\gsh$-deformation one can again associate a Kodaira-Spencer class, which now can be represented, essentially, by a pair $$(\phi,\psi)\in\g^1\otimes\m_S\oplus \h^0\otimes\m_S$$ where $\phi$ satisfies integrability as above, and the pair satisfy {\it{compatibility}} in the form $$\bp\psi=C(\ad(-\psi))(\phi)=
\sum\limits_{n=0}^\infty c_n\ad(-\psi)^n(\phi)$$ where the $c_n$ are the Bernoulli coefficients, defined by $$\frac{x}{e^x-1}=\sum\limits_{n=0}^\infty c_nx^n.$$ These conditions can again be formulated cohomologically via the deformation ring $R(\gsh)$, which is defined from the cohomology a suitable complex that we call {\it{Jacobi-Bernoulli }} complex, and which generalizes the Jacobi complex of a dgla.\par
This section is organized as follows. In \S2.1 we present basic definitions, properties and examples pertaining to deformations
of a Lie atom $\gsh$, viewed as torsors. The more familiar dgla case is also reviewed. In \S2.2 we define the Kodaira-Spencer class of a $\gsh$-deformation parametrized by $S$, and characterize it in terms of local homomorphisms $R(\gsh)\to S$, plus some extra data (the latter is not needed if $\gsh$ is a Lie pair).
 \subheading{2.1 Deformations as torsors}\par Our purpose here
is to define and study deformations with respect to a Lie atom
$\gsh=(\g,\h,i)$. Roughly speaking a $\gsh-$deformation consists of
a  $\g-$ deformation $\phi$, plus a 'trivialization of $\phi$ when
viewed as $\h^+-$ deformation'.\par We recall first the notion of
$\g-$deformation. Let $\g$ be a sheaf of Lie algebras over a
topological space $X$, let $E$ be a $\g-$module and $S$ a
finite-dimensional local $\C-$algebra with (nilpotent) maximal ideal
$\m$ of exponent $e$. We set $S_i=S/\m^{i+1}, \m_i=\m/\m^{i+1}.$
Then $\g\otimes\m$ has a natural structure of nilpotent $S$-Lie
algebra, with bracket coming from the bracket on $\g$ and the
multiplication $\m\times\m\to\m$. Note that there is a sheaf of
groups $G_S$ given by
$$G_S=\exp (\g\otimes\m)\tag 2.1.1$$(where $\exp=\exp_S$ is a finite series by
nilpotence of $\m$), with multiplication given by the
Campbell-Hausdorff formula, where $\exp$, as a map to
$\U_S(\g\otimes\m)$, is injective because the formal $\log$ series
gives an inverse.

\comment Though not essential for our purposes, it may be noted that
$G_S$ coincides with the (multiplicative) subgroup sheaf of sections
congruent to 1 modulo the ideal $\U^+ (\g\otimes \m)$ generated by
$\g\otimes\m$ in the universal enveloping algebra
$$\U(\g\otimes \m):=\U_S(\g\otimes \m).$$
This is easy to prove by induction on the exponent of $S$: note that
if $I<S$ is an ideal with $\m .I=0$ then $\g\otimes I\subseteq
\g\otimes\m$ is a central ideal yielding a central subgroup
$G_I=1+\g\otimes I\subseteq G_S$ and a central ideal $\g\otimes
I\subseteq\U(\g\otimes\m)$, hence a central subgroup $1+\g\otimes I$
in the multiplicative group of $\U(\g\otimes\m).$

\endcomment
\par
\proclaim{Definition 2.1.1} (i) A {\it{ $\g-$deformation over $S$}}
is a $G_S$-torsor, i.e. sheaf $G^{\phi}$ of sets with $G_S$-action
that is simply transitive, locally over $X$. An equivalence of
deformations is a $G_S$-equivariant map $G^\phi \to G^{\phi'}$ over
$X$.\par (ii) Given the $\g$-module $E$, a{\it{ $\g-$deformation of
$E$ over $S$}} is a sheaf $E^{\phi}$ of $S$-modules, together with a
maximal atlas of trivialisations
$$\Phi_{\alpha}:E^{\phi}|_{U_{\alpha}}\overset{\sim}\to{\to}
 E|_{U_{\alpha}}\otimes S,$$
such that the transition maps
$$\Psi_{\alpha\beta}:=\Phi_{\beta}\circ\Phi_\alpha^{-1}\in
\bar{G}_S(U_{\alpha}\cap U_{\beta})$$ where $\bar{\g}\subset\gl(E)$
is the image of $\g$ and $\bar{G}_S$ is the corresponding group
sheaf. \endproclaim\remark{Apology} The notation $G^\phi$
for a deformation is somewhat misleading because a deformation does not 'come with' a $\phi$ though a suitable $\phi$ does give rise to a deformation and any deformation comes from a $\phi$. Hopefully, the context will make our intention clear in each use of this notation.\endremark \remark{Remarks} (i) Given an abstract
$\g$-deformation $G^\phi$ and a $\g$-module $E$, a corresponding
$\g$-deformation $E^\phi$ of $E$ can be defined, either as the
deformation with the same transition functions, or as
$E^\phi=(E\times_X G^\phi)/G_S.$ Conversely, given a
$\g$-deformation $E^\phi$, it defines a $\bar{\g}$-deformation,
either as the deformation with the same transition functions or as
the sheaf of maps $E\otimes S\to E^\phi$ (or $E^\phi\to E\otimes S$)
that are locally in
$\bar{G}_S$.\par 
(ii) The notion of homomorphism of $\g-$deformations is defined in
the obvious way (with local representatives that are $\g-$ and
$S-$linear). The isomorphism class of  a $\g-$deformation is given
by the class of $(\Psi_{\alpha\beta})$ in the nonabelian \v{C}ech
cohomology set $\check{H}^1(X,G_S)$. \par (iii) Note that a
$\g-$deformation determines a functor from the category of
$\g-$modules to that of $\g-$deformations of modules. This functor
is exact and linear (in the sense that it induces a $\C-$linear map
from Hom$_\g(E_1, E_2)$ to Hom$_S(E_1^\phi,E_2^\phi)$. Moreover this
functor is determined by its value on any {\it{faithful}}
$\g-$module $E$ . We may call $E^\phi$ a model of $\phi$ or
$(\Psi_{\alpha\beta})$.\par (iv) Via the adjoint representation, any
$\g$-deformation $\phi$ over $S$ determines a $\g$-deformation of
$\g$ itself, $\g^\phi$, which is easily seen to be an $S$-Lie
algebra ($\g^\phi$ also coincides with the algebra of vertical
vector fields of $G^\phi/X$). In general, $\g^\phi$ need not
determine $\phi$.\par (v) If $(\g^., \bp)$ is a \dgla, one can
define an 'operatorial' notion of deformation as follows. Let
$\g^._{\bp}$ be the graded Lie algebra $\g^.\oplus\C[-1]\bp$, split
extension of $\C[-1]$ by $\g^.$ with bracket $$[\bp,a]=\bp(a),
[\bp,\bp]=0.$$ Then an operatorial $\g^.$-deformation is simply a
square-zero element of $\g^._{\bp}\otimes\m_S$ congruent to
$\bp\mod\m_S$, and two such are equivalent if they are in the same
orbit under the Adjoint action of $\exp(\g^0)$. The fact that, if
$\g^.$ is a suitable dgla resolution of $\g$, then equivalence
classes of $\g$-deformations coincide with equivalence classes of
operatorial $\g^.$-deformations is essentially the content of
Kodaira-Spencer theory (cf. \S2.2 below). \qed\endremark We now turn
to the case of Lie atoms and their deformations. Thus let
$\gsh=(\g,\h,i)$ be a sheaf of $S$-Lie atoms on $X$,
and let $E^\sharp =(E_1,E_2,j)$ be a sheaf of left $\gsh-$modules.
Note that $\gsh\otimes\m=(\g\otimes\m, \h\otimes\m,i\otimes\id)$ is
naturally a sheaf of $S$-Lie atoms, and we choose for its hull
$$\h^+_\m=({\text{\rom Lie\ closure\ of\ }}
\h\otimes\m) \subset\h^+\otimes\m.\tag 2.1.2$$ Note that the
inclusion $\h^+_\m\subset\h^+\otimes\m$ may well be strict, e.g. if
$\m^2=0$ then $\h^+_\m=\h\otimes\m$ always. Then $\h^+_\m$ is
endowed with the adjoint filtration as in (1.1.2):
$$\h\sbr 1._\m=\h\otimes\m\subset...\subset\h\sbr
i._\m=(\h\otimes\m)\sbr
i.\subset...\subset\h\sbr\infty._\m\subset\h^+_\m\tag 2.1.3$$

We have
 $$\h\sbr i._\m\otimes S_i=\h\sbr i+1._\m\otimes S_i=...
 =\h\sbr\infty._\m\otimes S_i .$$

 Working
with $\h^+_\m$ and $\h\sbr i._\m$, which depend on $\h$, rather than
$\h^+\otimes\m$, is the main purpose of specifying the
$\g$-submodule $\h\subset\h^+$ rather than working solely with
$\h^+.$ This is the main reason for introducing Lie atoms, as
opposed to working just with Lie pairs.\par Now set formally
$$H_S^+=\exp(\h^+_\m)\tag 2.1.4$$ (again, the exponential series is finite by
nilpotence of $\m$). By the Campbell-Hausdorff formula, $H^+_S$ is a
group, and there is a group homomorphism $G_S\to H^+_S.$ Moreover,
as noted in Lemma 1.1.3, the subset
$$H_S\sbr\infty.=\exp(\h_\m\sbr\infty.)\subset H^+_S\tag 2.1.5$$
is (left and right) $G_S$-invariant.  This notation should be used
with care because in general
$$H^+_S\subsetneqq (H^+)_S=\exp(\h^+\otimes \m_S).$$ Similarly,
$$H^+_{S_i}=\exp(\h^+_\m \otimes S_i).$$ Note that we have
 Lie pairs
$$\g^+_\m=(\g\otimes\m,\h^+_\m),$$
$$\g\sbr i._\m=(\g\otimes\m_i,\h^+_\m\otimes S_i)$$
and $\g\sbr i._\m=\g^+_\m$ for $i\geq e$ where $e$ is large enough
so $\m^{e+1}=0.$ Also set $$\gsh_\m=(\g\otimes\m,
\h_\m\sbr\infty.)=(\g\otimes\m, \h_\m\sbr e.).\tag 2.1.6$$ Thus
$$\gsh_{\m_i}=(\g\otimes\m_i, \h\sbr i._\m.)$$

\comment To an element $v\in\h\otimes\m$ we associate a formal
exponential
$$\exp(v)=\sum\limits_{i=0}^\infty\frac{(v)^i}{i!}
 ;$$ This may be viewed as an element of $\U(\h\otimes\m)$
 where $\h$ is viewed as an abelian Lie algebra. The set of these
 elements- which is just a copy of $\h\otimes \m$-
 is denoted $H_S.$
\endcomment
There is an obvious action
$$H^+_{S_i}\times E_2\otimes S_i \to E_2\otimes S_i,$$
$$\langle\exp(v),x\rangle=\sum\limits_{j=0}^i\frac{\ad(v)^j}{j!}
 (x)\tag 2.1.7$$
where $\ad(v)(x)=\langle v,x\rangle$.
Such maps are called {\it{left
$\h$-maps}} of $E^\sharp$ if $v\in\h_\m\sbr\infty.$. 
\comment Locally, an $S-$linear map
$$A:E_2\otimes S\to E_2\otimes S$$
is said to be a {\it{left $\h-$map}} if it is of the form
$$A=\exp(u)\circ A_v, u\in\g\otimes \m, v\in \h\otimes\m,$$
and similarly for right modules and right $\h-$maps. We note that
the set of left (resp. right) $\h-$maps is invariant under the left
(resp. right) action of $G_S$ on $\hom (E_1\otimes S,E_2\otimes S)$.
\endcomment
We consider the data of an $\h-$map to include
the element $v$ (this
is of course redundant if the action of $H_S$ is faithful), and two
such maps are considered equivalent if they belong to the same
$G_S-$orbit. Thus a left $\h-$ map is essentially a $G_S-$orbit of
an element of $\h\otimes\m$ in $H_S.$ Since $\m$ is nilpotent, any
left $\h-$map is an $S-$isomorphism. All the above leftist
considerations have obvious rightist analogues.

\par
The notion of $\h-$map
globalizes as follows. Given a  
 $\g-$deformation $E^\phi$,
a (global) {\it{left $\h-$map}} (with respect to $\phi$) is a map
$$A:E_2\otimes S\to E_2^{\phi}$$
such that for any atlas $\Phi_\alpha$ for $E_2^{\phi}$ over an open
covering $U_\alpha$, $\Phi_\alpha\circ A$
is given over $U_\alpha$
by a left $\h-$map. Note that this condition is independent of the
choice of atlas, and is moreover equivalent to the {\it existence}
of some atlas for which the $\Phi_\alpha\circ A$ are given by
$$x\mapsto \langle\exp(v_\alpha),x\rangle,
\ \  v_\alpha\in\h (U_\alpha)\otimes\m. \tag 2.1.8$$
We call such an
atlas a {\it{good atlas}} for $A$. Similarly, an 'abstract' or
'torsor' left $\h$-map is a map  $$A:(H^+)_S\to (H^+)^\phi\tag
2.1.9$$ from the trivial $H^+$ torsor to the $H^+$ torsor
determined by $\phi$, that is locally given by a left $\h$-map as
above. An $\h$-map as in (2.19) is equivalent to another such,
$$A':(H^+)_S\to(H^+)^{\phi'}$$ if there exists an equivalence of
$\g$-deformations $\epsilon:G^\phi\to G^{\phi'}$ and an element
$\gamma\in G_S$ such that, if we denote by $\epsilon^{H^+}$ the
natural extension of $\epsilon$ to an equivalence of
$\h^+$-deformations, then the following commutes
$$\matrix
(H^+)_S&\overset{A}\to{\to}&(H^+)^\phi\\
\gamma\downarrow&&\downarrow\epsilon^{H^+}\\
(H^+)_S&\overset{A'}\to{\to}&(H^+)^{\phi'}\endmatrix$$

\par The notion
of global right $\h-$map
$$B:F_1^\phi\to F_1\otimes S$$
for a right $\gsh-$module $(F_1, F_2, k)$ is defined similarly, as
is that of abstract global  $\h-$maps without
specifying a module. A
pair $(A,B)$ consisting of a left and right $\h-$map
is said to be a
{\it{dual pair}} if there exists a common
good atlas with respect to
which $A$ has the form (2.2) while $B$
has the form
$$x\mapsto   \langle\exp(-v_\alpha),x\rangle $$
with the {\it{same}} $v_\alpha$.
\proclaim{Definition 2.1.2} In the
above
situation, a left $\gsh-$deformation 
over $S$ consists of a $\g-$deformation $G^\phi$ together with a
left $\h-$map  from the trivial deformation to the
$\h^+$-deformation corresponding to $\phi$:
$$A:(H^+)_S\to (H^+)^{\phi}.$$
Similarly for right $\gsh-$deformation. A (2-sided)
$\gsh-$deformation consists of a $\g-$deformation $\phi$ together
with a dual pair $(A,B)$ of $\h-$maps with respect to $\phi$.
\endproclaim An obvious, yet fundamental observation is that the
various notions of $\gsh$-deformations are {\it functorial} with
respect to homomorphisms of Lie atoms. In particular, given a Lie
algebra $\k$ and a Lie homomorphism $\k\to\gsh$, any
$\k$-deformation over $S$ induces a $\gsh$-deformation over $S$.
\remark{Examples 1.1.4 cont.}
\item
{\bf A-B} When $E_1<E_2$ are vector spaces, $\gl(E_1<E_2)$-
deformation theory is just the local geometry of the Grassmannian
$G(\dim E_1, \dim E_2)$. When $E_1<E_2$ are vector bundles,
$\gl(E_1<E_2)$- deformation theory is just the local
geometry of the
'Grassmannian' of subbundles of $E_2$ or equivalently, the Quot
scheme of $E_2$ localized at the quotient $E_2\to E_2/E_1$.
Generally, if $A\subset B$ are coherent sheaves on a projective
scheme, we have the differential graded Lie atoms $\gl(A<B)$ and
$\gl(B>B/A)$, and for both of them the associated
deformation theory
is that of the Quot scheme of $B$ localized at $B/A$.\par
\item {\bf
C} When $\gsh = (\D^1(E),\D^2(E))$ is the heat algebra of the
invertible sheaf $E$, $\gsh-$ deformations of $E^\sharp=(E,E)$ are
called {\it{heat deformations}}. Recall that a $\D^1(E)-$
deformation consists of a deformation $\O^\phi$ of the structure
sheaf of $X$, together with an invertible
$\O^\phi-$ module $E^\phi$
that is a deformation of $E$. Lifting this to a
$\gsh-$ deformation
amounts to constructing $S-$linear, globally defined $\h$-maps
({\it{heat operators}})
$$A:E\otimes S\to E^\phi,$$
$$B:E^\phi\to E\otimes S.$$ By definition, $\h^+_\m$ is the
subalgebra of $\D^\infty(E)\otimes\m$ generated by
$\D^2(E)\otimes\m$, so any element of  $\h^+_\m$ is an $S$-linear
differential operator $\square$ that has order
$\leq i \mod\m^i$ for
any $i\geq 2.$ and $A$ and $B$ are locally
(with respect to an atlas
and a trivialisation of $E$) of the form
$$f\mapsto
\sum \frac{\square^i}{i!}f .$$ Note that this operator is of order
$\leq 2i-2 \mod \m^i, i\geq 2.$ Indeed writing locally
$\square=\sum\square_j, \square_j\equiv 0\mod\m^j,
\ord(\square_j)=j,$ we see that
$\square^i\equiv 0\mod\m^i$ and that
the highest-order term in $\square^{i-1}$ comes from
$\square_2^{i-1}$.
Moreover, $A,B$ yield mutally inverse $S$-linear
(not $\O_X$-linear)
isomorphisms. Notice that the heat operator $A$ yields a
well-defined lifting of sections (as well as cohomology classes,
etc.) of $E$ defined in any open set $U$ of $X$ to sections of
$E^\phi$ in $U$ (viz. $s\mapsto A(s\otimes 1)$). In particular,
suppose that $X$ is a compact
complex manifold 
Then
$$H^*(A):H^*(X,E)\otimes S\to H^*(X,E^\phi)$$
is an $S$-linear isomorphism.  Thus for any heat deformation the
cohomology module  $H^*(E^\phi)$ is not just (relatively)
unobstructed but {\it{canonically trivialised}}. Put another way,
$H^*(E^\phi)$ is endowed with a {\it{canonical flat connection}}
$$\nabla^\phi:H^*(E^\phi)\to H^0(E^\phi)\otimes\Omega_S\tag 2.8$$
determined by the requirement that
$$\nabla^\phi\circ H^0(A)(H^0(X,E)) =0,\tag 2.9$$
i.e. that the heat lift of sections of $E$ be flat.
\item {\bf D.} Here $\gsh = N_{Y/X}[-1]=(T_{X/Y},T_X).$ To start
with, a $T_{X/Y}-$ deformation is simply a deformation
in the usual
sense of the pair $(Y,X)$, giving rise, e.g. to a deformation
$(\I_Y^\phi,\O_X^\phi)$ of the left $\gsh-$module
$(\I_Y, \O_X)$ and
to a deformation $(\O_X^\phi, \O_Y^\phi)$ of the
right $\gsh-$module
$(\O_X,\O_Y)$. Then a left (resp. right) $\gsh-$ deformation of
$(\I_Y,\O_X)$ (resp. $(\O_X,\O_Y)$ consists of a
$T_{X/Y}-$deformation,  together with a $T_X-$map
$$A:\O_X^\phi\to \O_X\otimes S,$$
(resp. $B:\O_X\otimes S\to\O_X^\phi$).\nl Either $A$ or $B$ yield
trivialisations   of the deformation $\O_X^\phi$.
Thus left $\gsh-$
deformations yield deformations of $Y$ in a fixed
$X$, and similarly
for right deformations. Conversely, given a
deformation of $Y$ in a
fixed $X$, let $(x_\alpha^k)$ be local equations for $Y$ in $X$,
part of a local coordinate system. Then it is easy to see that we
can write equations for the deformation of $Y$ in the form
$$\exp(v_\alpha) (x_\alpha^k),\ \ v_\alpha\in T_X\otimes\m$$
($v_\alpha$ independent of $k$), so this comes from a left and a
right $\gsh-$ deformation of the form
$((\Psi_{\alpha\beta}),(v_\alpha))$ where
$$\Psi_{\alpha\beta}=\exp(v_\alpha)\exp(-v_\beta)\in
\U_S(T_{X/Y}\otimes\m) (U_{\alpha}\cap U_\beta).$$ Thus the four
notions of left, right and 2-sided $N_{Y/X}[-1]$-deformations and
deformations of $Y$ in a fixed $X$ all coincide.\nl Similarly, for a
mapping $f:Y\to X$, an $N_{f/X}[-1]$-
deformation is a deformation of
$(f,Y)$ fixing $X$. For $f$ submersive, a
$T_{f/Y}[-1]$-deformation
is a deformation of $(Y,f)$ fixing $Y$.
\item
{\bf E.} In this case we see similarly that $T_X\otimes\O_Y-$
deformations of $(\O_X\to\O_Y)$
consist of a deformation of
the pair $(X,Y)$, together with
trivialisations of the corresponding
deformations of $X$ and $Y$ separately, i.e. these are just
deformations of the embedding
$Y\hookrightarrow X$, fixing both $X$
and $Y$.

\endremark

\subheading{2.2 Kodaira-Spencer class}\par Our purpose here is to
associate a (higher-order) Kodaira-Spencer class to a
$\gsh$-deformation. In the next section we will show that such
classes can, under suitable hypotheses, be used to classify
$\gsh$-deformations.\par The basic idea of Kodaira-Spencer theory
can be described thus: let $(\gt^.,\bp)$ be a suitable acyclic
resolution of $\g$, let $G^0_S=\exp(\g^0\otimes\m)$. Then a
$\g$-deformation, i.e. a $G_S$-torsor, may be considered as a
subsheaf of the trivial $G^0_S$-torsor. This subsheaf can be
determined by specifying its tangent space, which can be determined
by a suitable operator on $\gt^.\otimes\m$ deforming $\bp.$ The
deformation of the operator is essentially the Kodaira-Spencer class
$\phi$ of the $\g$-deformation. Extending this to $\gsh$-deformation
is a matter of accounting in terms of $\phi$ for a trivialization of
the $\h^+$-deformation corresponding to $\phi$ . This is trickier
than one might expect, in part because the compatibility relation
characterizing an $\h^+$-trivialization of $\phi$ is not a linear or even quadratic condition (as is e.g. the integrability condition on
$\phi$), but an $n$-th degree condition, where $n$ is the order of
the deformation.\par Consider then a $\gsh$-deformation as in
Definition 2.1.2, consisting of a $G_S$-torsor $G^\phi$, with
transitions
$$\Psi_{\alpha\beta}=\exp (u_{\alpha\beta}),
u_{\alpha\beta}\in\g\otimes\m $$ with respect to a suitable open
covering $(U_\alpha)$,  plus an $\h$-map
$$A:(H^+)_S\to (H^+)^\phi$$
locally given as $$\exp(v_\alpha),
v_{\alpha\beta}\in\h_\m\sbr\infty..$$
Then the condition that the $\exp(v_\alpha)$
should glue together to
a globally defined map left $\h-$ map $A$ is
$$\Psi_{\alpha\beta}\circ \exp(v_\alpha )=
\exp(v_\beta).\tag 2.2.1$$
By analogy with the procedure of \cite{R, p.61},
this condition may
be analyzed in terms of Kodaira-Spencer cochains, as follows.Let
$$\tilde{\gsh}=(\gt^.,\tilde{\h}^., \bp)$$
be a suitable acyclic (soft or injective) \dgla resolution of $\g$
and set
$$(\g^.)^\sharp=(\g^.,\h^.,\bp):=
(\Gamma(\gt^.),\Gamma(\tilde{\h}^.),
\bp)$$ Write
$$\Psi_{ab}=\exp(-s_\alpha)\exp(s_\beta),$$
$$s_\alpha\in\gt^0(U_\alpha)\otimes\m,$$ where the cochain
$(\exp(s_\alpha))$ is determined up to left multiplication by an
element of the 'soft gauge group' $$\exp(\mu)\in
G^0_S=\exp(\Gamma(\gt^0)\otimes\m), $$ with $\mu$ independent of
$\alpha$.
 Thus
$$\exp(-s_\alpha)\exp(v_\alpha)=\exp(-s_\beta)\exp(v_\beta)$$ so
this element is globally defined and, by the Campbell-Hausdorff
formula, can be written as
$$\exp(\psi),\ \  \psi\in(\h^0_\m)\sbr\infty.=(\h^0_\m)\sbr e.$$
(where $e$ is such that $\m^{e+1}=0$).
As observed in \cite{Rcid}, the element
$$\phi:=\exp(-s_\alpha)\bp\exp(s_\alpha)\in\g^1\otimes\m$$
is globally
defined independent of $\alpha$ and is the
Kodaira-Spencer cochain
defining the deformation, and we can write formally
$$\phi=D(-\ad(s))(\bp s)$$ where
$$D(x)=(e^x-1)/x$$ is the Deligne function.
$\phi$ is well-defined up to replacing $\exp(s_\alpha)$ by
$\exp(\mu)\exp(s_\alpha)$ as above (
$\mu$ independent of $\alpha$),
which is equivalent to conjugating the operator $\bp+\phi$ by
$\exp(\mu)$. Consider the Bernoulli function
$$C(x)=1/D(x)=\sum\limits_{n=0}^\infty c_nx^n.$$
Thus $c_n=B_n/n!$ where
$B_n$ is the $n$th Bernoulli number
(cf. \S0.2) . \comment Then it
is elementary that $c_0=1$ and
$$c_n=\sum\limits_{1\leq i\leq m\leq
n}(-1)^i\binom{m}{i}\frac{i^n}{(m+1)n!},
\forall n\geq 1.\tag ??$$ In
particular, $c_{2n+1}=0,\forall n\geq 1.$  $\lceil$ Indeed set
$y=e^x-1$ so
$C(x)=\log(1+y)/y=\sum\limits_{m=0}^\infty(-1)^my^m/(m+1)$,
and use
the fact that ord$_x(y^m)=m$, so that
$$y^m=\sum\limits_{i=0}^m\sum\limits_{n\geq
m}(-1)^{i-m}\binom{m}{i}\frac{(ix)^n}{n!}.\rfloor$$\par
\endcomment
 Then a formal calculation in
$\U_S(\Gamma(\h^{.})^+_\m)$ shows that
$$\bp\exp(\psi)=-\phi\exp(\psi),\bp\exp(-\psi)=\exp(-\psi)\phi$$
This implies
$$i(\phi)=\bp(\exp(\psi))\exp(-\psi)=
D(-\ad(\psi))(\bp\psi)$$ hence
$$\bp\psi=C(\ad(-\psi))i(\phi)=
\sum\limits_{n=0}^\infty c_n\ad(-\psi)^ni(\phi).\tag 2.2.2$$ In
other words, the vector $$(-\psi, \phi,
-\psi\otimes\phi,...,(-\psi)^n\otimes\phi,...)\in (\h^0_\m)\sbr
e.\oplus\g^1\otimes\m\oplus ...\oplus\Sym^n((\h^0_\m)\sbr
e.)\otimes\g^1\otimes\m$$ {\it is a cocycle
for the JacoBer complex}
$\Jsh(\gsh_\m)$. In view of the definition of the
JacoBer complex,
it follows that for any $i\geq 1$, the vector
$$\epsilon(\phi,\psi)=
(\bigwedge\limits^r\phi\otimes(-\psi)^n)\in
\bigoplus\limits_{r+n\leq
i} \bigwedge\limits^r(\g^1\otimes\m_i)\otimes
\Sym^n((\h^0_{\m_i})\sbr i.)),\tag 2.2.3$$ which is a priori a
0-cochain for the complex $\Jsh_i(\gsh)$,
is a cocycle as well. The
associated cohomology class $$\alpha_i=\alpha_i(\phi,\psi) \in
\H^0(\Jsh_i(\gsh_{\m_i}))\subset
\H^0(\Jsh_i(\gsh))\otimes \m_i\tag
2.2.4$$ is called the $i$-th {\it Kodaira-Spencer class}  of the
$\gsh$-deformation $(G^\phi, A)$. For all $i\geq e$, clearly
$\alpha_i$ has the same value, which we denote by
$\alpha(\phi,\psi)$ and call 'the' Kodaira-Spencer class of the
deformation. It is also clear from the definitions that the
$\alpha_i(\phi,\psi)$ are
comultiplicative or 'morphic', hence their
'transpose' yields a sequence of ring homomorphisms
$$^t\alpha_i(\phi,\psi):R_i(\gsh)=
\C\oplus\H^0(\Jsh_i(\gsh))^*\to S_i\tag 2.2.5$$ that are mutually
compatible via (1.2.6).  Obviously, $^t\alpha_i(\phi,\psi)$
determines $\alpha_i(\phi,\psi)$. However, it is not clear a priori
that given a homomorphism as in (2.2.5), it necessarily comes from
an element of $\H^0(\Jsh_i(\gsh_{\m_i}))$, as opposed to $
\H^0(\Jsh_i(\gsh))\otimes \m_i$.
\par Nonetheless, suppose conversely that we have a local artinian
$\C$-algebra $S$ of exponent $e$ and a compatible sequence of
homomorphisms
$$\beta_i:R_i(\gsh)\to S_i, i\leq e.$$
Suppose moreover, as in Definition 1.2.3, that $\gsh$ is positive.
Then we may further assume that $\g^{\leq 0}=0, \h^{<0}=0.$ Clearly
each $\beta_i$ can be written as as $^t\alpha_i(\phi_i, \psi_i)$
where
$$\alpha_i(\phi_i, \psi_i)\in\H^0(\Jsh_i(\gsh))\otimes \m_i.$$
Thanks to our vanishing  hypothesis, the pairs $\phi_i, \psi_i$ are
{\it uniquely determined}. Therefore
$$\psi_i\equiv \psi_e\mod\m^{i+1}.$$ Since $\psi_i\in(\h^0)\sbr
i.\otimes m_i$, this implies that $$\psi_e\in (\h^0_\m)\sbr e.=
(\h^0_\m)\sbr \infty.,$$ so that $$\alpha_e(\phi_e, \psi_e)\in
\H^0(\Jsh _e(\gsh_\m))
 =\H^0(\Jsh(\g\otimes\m,\h_\m \sbr\infty.)),$$
and clearly
$$\beta_i=^t\alpha_i(\phi_e, \psi_e).$$ We have established the
following\proclaim{Proposition 2.2.1} Let $\gsh=(\g,\h)$ be a
positive Lie atom,  and let $S$ be a local artinian $\C$-algebra of
exponent $e$. Then\par (i) there is a 1-1 correspondence between
compatible sequences of morphic elements
$$\alpha_i \in \H^0(\Jsh_i(\gsh_{\m_i})), i=1,...,e$$
and compatible
sequences of homomorphisms $$\beta_i:R_i(\gsh)\to S_i, i\leq e;$$
(ii) the Kodaira-Spencer class of any $\gsh$-deformation over $S$
yields such sequences.

\endproclaim We call a sequence $\alpha.=(\alpha_1,...,\alpha_e)$ or
$\beta.=(\beta_1,...,\beta_e)$ a {\it Kodaira-Spencer sequence} for
$\gsh$ and $S$. The following is obvious\proclaim{Corollary 2.2.2}
(i)  A Kodaira-Spencer sequence is uniquely determined by its last
element.
\par (ii) A local homomorphism
$\beta_\infty:R(\gsh)=R_\infty(\gsh)\to S$ gives rise to a
Kodaira-Spencer sequence $\beta.$ iff for each $i$, $\beta_\infty$
maps the kernel of $R(\gsh)\to R_i(\gsh)$, an ideal which a priori
contains $\m_{R(\gsh)}^{i+1}$, to $\m_S^{i+1}.$\par (iii) If $\gsh$
is a Lie pair, any local homomorphism $\beta_\infty:R(\gsh)\to S$
yields a Kodaira-Spencer sequence.
\endproclaim
\remark{Remark 2.2.3} When $\gsh$ is not a Lie pair, $R_i(\gsh)$
will in general be a proper quotient of $R_e(\gsh)/\m_e^{i+1}$, so
assertion (iii) above will not hold.
\endremark

\comment
 which is equivalent
to the following equation in $\U(\g\otimes\m,\h\otimes\m)$, in which
we set
$$D(x)=\frac{\exp (x)-1}{x}=
\sum\limits_{k=1}^\infty \frac{x^k}{(k+1)!}\ \  :\tag 2.4$$
$$D(u_{\alpha\beta})i(u_{\alpha\beta})+
\exp (u_{\alpha\beta}).v_\alpha =v_\beta . \tag 2.5$$ The condition
for a right $\h -$ map $B$ to be given by the $\id-v_\alpha$ is
$$i(u_{\alpha\beta})D(u_{\alpha\beta})+
v_\alpha\exp (u_{\alpha\beta}) =v_\beta . \tag 2.6$$

\remark{Examples 1.3 cont.}
\item{\bf A-B} When $E_1<E_2$ are vector spaces, $\gl(E_1<E_2)$-
deformation theory is just the local geometry of the Grassmannian
$G(\dim E_1, \dim E_2)$. When $E_1<E_2$ are vector bundles,
$\gl(E_1<E_2)$-
deformation theory is just the local geometry of the
'Grassmannian' of subbundles of $E_2$ or equivalently, the Quot
scheme of $E_2$ localized at the quotient $E_2\to E_2/E_1$.
Generally, if $A\subset B$ are coherent sheaves on a projective
scheme, we have the differential graded Lie atoms $\gl(A<B)$ and
$\gl(B>B/A)$, and for both of them the associated deformation theory
is that of the Quot scheme of $B$ localized at $B/A$.\par
\item {\bf
C} When $\gsh = (\D^1(E),\D^2(E))$ is the heat algebra of the
invertible sheaf $E$, $\gsh-$ deformations of $E^\sharp=(E,E)$ are
called {\it{heat deformations}}. Recall that a $\D^1(E)-$
deformation consists of a deformation $\O^\phi$ of the structure
sheaf of $X$, together with an invertible
$\O^\phi-$ module $E^\phi$
that is a deformation of $E$. Lifting this to a $\gsh-$ deformation
amounts to constructing $S-$linear, globally defined $\h$-maps
({\it{heat operators}})
$$A:E\otimes S\to E^\phi,$$
$$B:E^\phi\to E\otimes S.$$ By definition, $\h^+_\m$ is the
subalgebra of $\D^\infty(E)\otimes\m$ generated by
$\D^2(E)\otimes\m$, so any element of  $\h^+_\m$ is an $S$-linear
differential operator $\square$
that has order $\leq i \mod\m^i$ for
any $i\geq 2.$ and $A$ and $B$ are locally
(with respect to an atlas
and a trivialisation of $E$) of the form
$$f\mapsto
\sum \frac{\square^i}{i!}f .$$ Note that this operator is of order
$\leq 2i-2 \mod \m^i, i\geq 2.$ Indeed writing locally
$\square=\sum\square_j, \square_j\equiv 0\mod\m^j,
\ord(\square_j)=j,$ we see that $\square^i\equiv 0\mod\m^i$ and that
the highest-order term in $\square^{i-1}$ comes from
$\square_2^{i-1}$.
Moreover, $A,B$ yield mutally inverse $S$-linear (not $\O_X$-linear)
isomorphisms. Notice that the heat operator $A$ yields a
well-defined lifting of sections (as well as cohomology classes,
etc.) of $E$ defined in any open set $U$ of $X$ to sections of
$E^\phi$ in $U$ (viz. $s\mapsto A(s\otimes 1)$). In particular,
suppose that $X$ is a compact
complex manifold 
Then
$$H^*(A):H^*(X,E)\otimes S\to H^*(X,E^\phi)$$
is an $S$-linear isomorphism.  Thus for any heat deformation the
cohomology module  $H^*(E^\phi)$ is not just (relatively)
unobstructed but {\it{canonically trivialised}}. Put another way,
$H^*(E^\phi)$ is endowed with a {\it{canonical flat connection}}
$$\nabla^\phi:H^*(E^\phi)\to H^0(E^\phi)\otimes\Omega_S\tag 2.8$$
determined by the requirement that
$$\nabla^\phi\circ H^0(A)(H^0(X,E)) =0,\tag 2.9$$
i.e. that the heat lift of sections of $E$ be flat.
\item {\bf D.} Here $\gsh = N_{Y/X}[-1]=(T_{X/Y},T_X).$ To start
with, a $T_{X/Y}-$ deformation is simply a deformation in the usual
sense of the pair $(Y,X)$, giving rise, e.g. to a deformation
$(\I_Y^\phi,\O_X^\phi)$ of the left $\gsh-$module $(\I_Y, \O_X)$ and
to a deformation $(\O_X^\phi, \O_Y^\phi)$ of the right $\gsh-$module
$(\O_X,\O_Y)$. Then a left (resp. right) $\gsh-$ deformation of
$(\I_Y,\O_X)$ (resp. $(\O_X,\O_Y)$ consists of a
$T_{X/Y}-$deformation,  together with a $T_X-$map
$$A:\O_X^\phi\to \O_X\otimes S,$$
(resp. $B:\O_X\otimes S\to\O_X^\phi$).\nl Either $A$ or $B$ yield
trivialisations   of the deformation $\O_X^\phi$. Thus left $\gsh-$
deformations yield deformations of $Y$ in a fixed $X$, and similarly
for right deformations. Conversely, given a deformation of $Y$ in a
fixed $X$, let $(x_\alpha^k)$ be local equations for $Y$ in $X$,
part of a local coordinate system. Then it is easy to see that we
can write equations for the deformation of $Y$ in the form
$$\exp(v_\alpha) (x_\alpha^k),\ \ v_\alpha\in T_X\otimes\m$$
($v_\alpha$ independent of $k$), so this comes from a left and a
right $\gsh-$ deformation of the form
$((\Psi_{\alpha\beta}),(v_\alpha))$ where
$$\Psi_{\alpha\beta}=\exp(v_\alpha)\exp(-v_\beta)\in
\U_S(T_{X/Y}\otimes\m) (U_{\alpha}\cap U_\beta).$$ Thus the four
notions of left, right and 2-sided $N_{Y/X}[-1]$-deformations and
deformations of $Y$ in a fixed $X$ all coincide.\nl Similarly, for a
mapping $f:Y\to X$, an $N_{f/X}[-1]$-deformation is a deformation of
$(f,Y)$ fixing $X$. For $f$ submersive, a $T_{f/Y}[-1]$-deformation
is a deformation of $(Y,f)$ fixing $Y$.
\item
{\bf E.} In this case we see similarly that $T_X\otimes\O_Y-$
deformations of $(\I_Y\oplus\O_Y,\O_X)$ consist of a deformation of
the pair $(X,Y)$, together with trivialisations of the corresponding
deformations of $X$ and $Y$ separately, i.e. these are just
deformations of the embedding $Y\hookrightarrow X$, fixing both $X$
and $Y$.

\endremark
\endcomment
\heading 3. Universal deformations\endheading
A principal goal of $\g$-deformation theory is the construction of a universal $\g$-deformation, i.e. a $\g$-deformation over $R(\g)$,  from which any $\g$-deformation over $S$ is obtained via a (unique) base-change map $R(\g)\to S$. [Actually what one seeks is, for each $e\geq 0$, an $e$-universal deformation, which is one having the above universality property for all $S$ of exponent $e$, i.e. such that $\m_S^{e+1}=0$; then one can take the limit as $e\to\infty$ to get the (formally) universal deformation. But we will ignore such technicalities in this introductory paragraph.] Proving the existence of a universal deformation typically requires some restrictive hypotheses on $\g$, along the lines of nonexistence of automorphisms $(H^0(\g)=0)$, as well as some more technical finiteness hypotheses ('admissibility'). See \cite{Rcid, Ruvhs} for details. The proof typically proceeds according to the following outline:\par
(i) 'inverting' the Kodaira-Spencer class, e.g by associating to a local homomorphism $$h:R(\g)\to S$$ a $\g$-deformation $\beta(h)$ over $S$; in particular, by applying this to the identity map of $R(\g)$ we obtain a 'distinguished' deformation $\phi^u$ over $R(\g)$ (that one would like to prove is actually universal);\par (ii) proving that for any $\g$-deformation $\phi$ over $S$, the Kodaira-Spencer homomorphism $$\alpha(\phi):R(\g)\to S$$ constructed above depends only on the equivalence (i.e. $G$-conjugacy) class of $\phi$ as deformation;\par (iii) proving that for any local homomorphism $h$ as above, the Kodaira-Spencer homomorphism associated to $h^*(\phi^u)$ is just $h$;\par (iv)
 proving that for any $\g$-deformation $\phi$, $\phi$ is equivalent to $\alpha(\phi)^*(\phi^u)$.\par
 The conjunction of (i)-(iv) implies that the assignment
 $$h\mapsto h^*(\phi^u)$$
 yields a bijection
 $$\{{\text{local }}\C-{\text {algebra homomorphisms }} R(\g)\to S \}\to \{\g{\text{-deformations over }} S\}$$ which implies the universality of $\phi^u$.\par
The main purpose of this section is to
construct, under suitable hypotheses, the universal $\gsh$-deformation associated to a sheaf $\gsh$ of
Lie atoms, which is simultaneously the universal $\gsh-$deformation
of any $\gsh-$module $E^\sharp$ (see Theorem 3.3 below). This universal deformation is the evident analogue of the corresponding notion for $\g$-deformations. We thus extend the main result of
[Rcid] to the cases of atoms. We shall also prove a generalization (see Theorem 3.1) of the latter result where the hypothesis of trivial sections of weakened to 'central sections', i.e. that $H^0(\g)$ maps to the center of $\g^.(U)$ for all open sets $U$.
\par\subheading{3.1 dgla case revisited}
We shall assume throughout, without explicit mention, that all
sheaves of Lie algebras and modules considered are {\it admissible}
in the sense of [Rcid]. In addition, unless otherwise stated we
shall assume their cohomology is finite-dimensional. We begin by
reviewing the main construction of [Rcid] and restating its main
theorem in a stronger form. As in \S1, we have the
{\it{Jacobi complex}} $J_m(\g)$,
which is identified with the Jacobi
complex associated to the \dgla $\g^.=(\Gamma(\gt^.), \bp)$, where
$(\gt^., \bp)$ is a suitable acyclic (injective, flabby, soft)
resolution. This is
a complex in degrees $[-m,-1]$ 
\comment parametrizing nonempty subsets of $X$ of cardinality at
most $m$, which has the form
$$\lb^m(\g)\to ...\to\lb^2(\g)\to \g$$
where $\lb^k(\g)$ is the exterior alternating tensor power and the
coboundary maps
$$d_k:\lb^k(\g)\to\lb^{k-1}(\g)$$
are given by
$$d_k(a_1\times ...\times a_k)=
\frac{1}{k!}\sum\limits_{\pi\in S_k} [\sgn (\pi )
[a_{\pi (1)},a_{\pi
(2)}] \times a_{\pi (3)}...\times a_{\pi (k)}. \tag 3.1$$  We showed
in [Rcid ] that $J_m(\g)$ \endcomment
which admits a
comultiplicative structure, making
$$R_m(\g)=\C\oplus\H^0(J_m(\g))^*$$
 a $\C-$algebra (finite-dimensional by the admissibility
hypothesis), and we constructed a certain 'tautological'
$\g-$deformation $u_m$ over it (more precisely, $u_m$ is only
defined up to equivalence- more details below). To any
$\g-$deformation $\phi$ over an algebra $(S,\m)$ of exponent $m$
 we associated a canonical Kodaira-Spencer
homomorphism
$$\alpha =\alpha (\phi ) :R_m(\g )\to S; \ $$ e.g. $u_m$ is
essentially characterized by the property that
$\alpha(u_m)=\id_{R_m(\g)}.$ Although in [Rcid] we made the
hypothesis that $H^0(\g)=0,$ this is in fact not needed for the
foregoing statements, and is only used in the proof that $u_m$ is an
$m-$universal deformation.\par Now the hypothesis $H^0(\g)=0$ can be
relaxed somewhat. Let us say that $\g$ has {\it{central sections}}
if for each open set $U\subset X$, the image of the restriction map
$$H^0(\g)\to \g(U)$$
is contained in the center of $\g(U)$. Equivalently, in terms of a
dgla resolution as above $\g\to\gt^.,$ the condition is that
$H^0(\g)$ be contained in the center of $\Gamma (\gt^.)$, i.e. the
bracket
$$H^0(\g )\times \Gamma (\gt^i)\to\Gamma (\gt^i)\tag 3.1$$
should vanish. \proclaim{Theorem 3.1} Let $\g$ be an admissible dgla
and suppose that $\g$ has central sections. Then for any
$\g-$deformation $\phi$ there exists an equivalence of deformations
$$\phi\simto\alpha (\phi )^*(u_m)=
u_m\otimes_{R_m(\g)}S;\tag 3.2$$ any two such equivalences
differ by
an element of
$${\text{ Aut}} (\phi )=H^0(\exp (g^\phi\otimes\m)).$$
In particular, if $H^0(\g)=0$ then the equivalence is unique.
Consequently, for any admissible pair $(\g,E)$ there are
equivalences
$$E^\phi\simto\alpha (\phi )^*(E^{u_m})$$
any two of which differ by an element of ${\text{ Aut}} (\phi ).$

\endproclaim
\demo{proof} We first prove the isomorphism (3.2). More generally,
we will show \proclaim{Lemma 3.2} For any two deformations $\phi_1,
\phi_2$, if $\epsilon(\phi_1)$ and $\epsilon(\phi_2)$ are
cohomologous, then $\phi_1, \phi_2$ are equivalent as
deformations.\endproclaim This Lemma yields the isomorphism (3.2) as
a special case because $\epsilon(\phi)$ and\nl $\epsilon(\alpha
(\phi )^*(u_m))$ are clearly cohomologous. The Lemma is a
generalization of \cite{Rcid}, Theorem 0.1, Step 4, pp. 63-64. We
will give a new, self-contained proof. This
 proof is by induction on
the exponent $e$ of $S$. Using induction, we may assume $\phi_1,
\phi_2$ are equivalent $\mod \m^e$. Hence there exist deformations
$\phi'_1, \phi'_2\in\g^1\otimes\m$ equivalent, respectively, to
$\phi_1, \phi_2$, such that $$\phi'_1\equiv\phi'_2\mod\m^e.$$
Because equivalent deformations yield cohomologous Kodaira-Spencer
classes (Step 2, p.62 of \cite{Rcid})  it follows that
$$[\epsilon(\phi'_1)]=[\epsilon(\phi_1)]=
[\epsilon(\phi_2)]=[\epsilon(\phi'_2)].$$ Replacing
$\phi_1, \phi_2$
by $\phi'_1,\phi'_2$, we may in fact assume that
$$\phi_1\equiv\phi_2\mod\m^e.$$ This implies that $$\phi_1^{
i}=\phi_2^{ i}\in\Sym^i(\g^1)\otimes\m^i, \forall i>1.$$ Our
assumption that $\epsilon(\phi_1)$ and $\epsilon(\phi_2)$ are
cohomologous implies that there exist $$\mu_i\in
\g^0\otimes\Sym^{i-1}(\g^1)\otimes\m, i=0,...,m,$$ such that
$\epsilon(\phi_1)-\epsilon(\phi_2)$ is the coboundary of $(\mu.)$.
Recall that the total coboundary of $J$ consists of a vertical part,
induced by $\bp$, and a horizontal part induced by the bracket. We
can write
$$\mu_i=\mu'_i\otimes\mu_i^", \forall i>0$$ with
$\mu_i^"\in\Sym^{i-1}(\g^1)\otimes\m$ linearly independent.
Working
backwards from the degree-$m$ component,
the fact that the part of
the coboundary of $(\mu.)$ in $\Sym^m(\Gamma(\g^1))$
is zero implies
that the vertical, ($\bp$-induced) coboundary of $\mu_m$ vanishes,
i.e.
$$\bp(\mu'_m)=0.$$ By Central Sections, it follows that the
horizontal (bracket-induced) coboundary of $\mu_m$ vanishes.
Therefore, the vertical coboundary of $\mu_{m-1}$ vanishes, etc.
Continuing backwards in this manner,
we conclude eventually that the
horizontal coboundary of $\mu_1$ is zero and
$$\bp(\mu_0)=\phi_1-\phi_2.$$ In particular, $\mu_0\equiv
0\mod\m^e$, so that $[\mu_0,\phi_1]=[\mu_0,\phi_2]=0.$ Therefore
clearly $\phi_1$ and $\phi_2$ are equivalent as deformations, as
$$\exp(-\mu_0))(\bp+\phi_2)\exp(\mu_0)=\bp+\phi_1.$$
QED Lemma 3.2.\par \comment

**********************
 This is a matter of adapting the argument in
the proof of Theorem 0.1, Step 4, pp.63-64 in [Rcid], and we will
just indicate the changes. We work in $\H^0(J_m(\g)) \otimes\m$
rather than $\H^0(J_m(\g), \m^.)$, which may not inject to it. Then,
with the notation of {\it{loc. cit.}} we may write
$$u_1=\sum v_i\otimes\phi_i\in \Gamma (\g^0)\otimes
\Gamma (\g^1)\otimes\m .$$ The argument there given shows that
$$u_1=u_0\times \phi + w\times\phi$$
where $w\times\phi\in H^0(\g)\otimes\Gamma (\g^1)\otimes \m.$ Now-
and this is the point- since
$$[w,\phi]=0$$
thanks to our assumption of central sections,
this is sufficient to
show that $\tilde{\phi}-\phi$ is the total coboundary of
$u_0\times
\phi$, as required.
\par
\endcomment
 From Lemma 3.2 we deduce the existence of an isomorphism
as in (3.2). Given this, the fact that two such isomorphisms
differ
by an element of ${\text{ Aut}} (\phi )$ is obvious.
To identify the
latter group it suffices to identify its Lie algebra
${\text {\rom
aut}}(\phi)$, which is given locally by the set of
$\g-$endomorphisms
$$\ad (v)\in\g^0\otimes\m$$ of the
resolution
$$(\g^.\otimes S, \bp +\ad (\phi )).$$
It is elementary to check that the condition on $v$ is precisely
$$\bp (v) +\ad (\phi )(v) =0,$$
i.e. $v\in \g^\phi$. Thus the local endomorphism algebra of $\phi$
is $\g^\phi$ and the global one is ${\text{\rom
aut}}(\phi)=H^0(\g^\phi)$ Finally, note that $\g^\phi$ admits a
Jordan-H\"older series with each subquotient isomorphic as a sheaf
to $\g$. Therefore if $H^0(\g)=0,$, we have
$$H^0(\g^\phi)=0, $$
hence ${\text{ Aut}} (\phi )=(1)$. \qed

\enddemo
\remark{Remark} Without the hypothesis of central sections it is
still possible to 'classify' $\g-$deformations over $(S,\m)$ in
terms of $\H^0(J_m(\g), \m^.)$ but it is not immediately clear how
this is related to semiuniversal deformations.
\endremark
\subheading{3.2 Case of Lie atoms}
 We now extend these results from Lie algebras
to Lie atoms. \proclaim{Theorem 3.3} Let $\gsh=(\g,\h)$ be an
admissible Lie atom such that $\g$ has central sections.
\par (i) For any $\gsh$-deformation $(\phi,\psi)$ over an artin
local $\C$-algebra $S$, the associated Kodaira-Spencer sequence
$\beta.(\phi,\psi)$ depends only on the equivalence class of
$(\phi,\psi)$; conversely, $\beta.(\phi,\psi)$ determines the
equivalence class of $(\phi,\psi)$.\par (ii) Assume $\gsh$ is a Lie
pair. Then for any natural number $m$ there exists an equivalence
class
$$u_m=(\phi^u_m,\psi^u_m)$$ of $\gsh$ deformations over
$R_m(\gsh)$ such that for
any local artin $\C$-algebra of exponent $\leq m$, the assignments
$$(\phi,\psi)\mapsto \beta(\phi,\psi),$$
$$\beta\mapsto \beta^*(u_m)$$ establish a bijection between the set
of equivalence classes of $\gsh$-deformation over $S$ and \nl${\text
{\rom{Hom}}}_{\C-\text{\rom alg}}(R_m(\gsh), S).$
\endproclaim
  \demo{proof} (i) The proof of these assertions is the same as
 in the  dgla case, Theorem 3.1 above.
 \par
(ii) Set $S_m=R_m(\gsh)$ with its maximal ideal $\m_m.$ Then the
identity map of $S_m$ corresponds to a morphic element of
$\H^0(\Jsh_m(\gsh\otimes\m_m))$, which can be written in the form
$\alpha_m(\phi^u_m,\psi^u_m)$, and the fact that
$\alpha_m(\phi^u_m,
\psi^u_m)$ comes from a cocycle implies that
$$u_m=(\phi_m^u,
\psi_m^u)\in(\Gamma(\g^1)\oplus\Gamma(\h^0)\otimes\m_m$$ is a
$\gsh$-deformation, which we call a {\it universal $m$-th order
$\gsh$ deformation}. By (i), the equivalence class of $u_m$ is
independent of the choice of representative for the identity.\par
Now given an arbitrary $\gsh$-deformation $(\phi,\psi)$ over an
artin local $\C$-algebra $(S,\m)$ of exponent $e$, we get as in
\S2.2 a Kodaira-Spencer homomorphism
$$\alpha=\alpha_e(\phi,\psi):R_e(\gsh)\to S.$$ The proof that
$\alpha$ depends only on the equivalence class of $(\phi,\psi)$ is
as in the proof of Theorem 3.1, as is the proof that
$(\phi, \psi)$
is equivalent to $\alpha^*(u_e).$ \qed
\enddemo
\remark{Remark} For a general, say positive, Lie atom $\gsh$,
 we are
not asserting the existence of a universal $\gsh$-deformation. If
$\g^+$ denotes the hull of $\gsh$, $\gsh$ deformations over $S$ are
classified by the set of maps $R_e(\g^+)\to S$ that happen to be compatible
with a Kodaira-Spencer sequence. But it's not clear that there is a
single algebra $R$ such that $\gsh$-deformations over
$S$ correspond
bijectively with ${\text{\rom {Hom}}}_{\C-\text{\rom{alg}}}(R,S).$
\endremark
 \comment ****************As observed in \S1, a Lie
atom $\gsh$ admits a Jacobi complex $J_m(\gsh)$ which, as in the
case of Lie algebras, admits a comultiplicative structure, which
amounts to the map

\comment
 First
recall the modular Jacobi complex $J_m(\g,\h)$ associated to a
$\g-$module $\h$. This is a complex in degrees $[-m,0]$ defined on a
space $X<m,1>$ parametrizing pointed subsets of $X$ of cardinality
at most $m+1.$ It has the form
$$\lb^m(\g)\bt\h\to ...\to\g\bt\h\to \h$$
with differentials
$$\partial_k(a_1\times ...\times a_k\times v)=
d_k(a_1\times ...\times a_k) \times v
-\frac{2}{k}\sum\limits_{j=i}^k (-1)^ja_1 \times ...\hat{a_j}\times
...\times a_k\times \langle a_j,v\rangle),$$ where $d_k$ is the
differential in $J_m(\g)$ (see [Ruvhs]).\par Now let
$\gsh=(\g,\h,i)$ be a Lie atom. Then the $\g-$homomorphism $i$ gives
rise to a map of complexes
$$i_m:J_m(\g)\to \pi_{m-1,1\ *}J_{m-1}(\g,\h),$$
where $ \pi_{m-1,1}:X<m-1,1>\to X<m>$ is the natural map. Then for
any $j>0$ the map $i$ induces $\g\bt\sigma^j\h\to\sigma^{j+1}\h$,
whence a map as in (1.2) $$\pi_{m-j,j\ *}J_m(\g,\sigma^j\h)[j-1]\to
\pi_{m-j-1,j+1\ *}J_m(\g, \sigma^{j+1}\h)[j]$$ Assembling these
together, we get $J_m(\gsh)$, the Jacobi complex of the atom $\gsh$.
We note the natural map

$$\sigma_m :
J_m(\gsh)\to (J_m/J_1)(\gsh)\to\Sym^2(J_{m-1}(\gsh))$$ obtained by
assembling together various 'exterior comultiplication' maps
$$\lb^j(\gsh)\to \lb^r(\gsh)\bt\lb^{j-r}(\gsh).$$

This gives rise to a  comultiplicative structure on
$H^0(J_m(\gsh))$, whence a local finite-dimensional $\C-$algebra
structure on
$$R_m(\gsh):= \C\oplus H^0(J_m(\gsh))^*$$
as well as a (multiplicative) local homomorphism
$$R_m(\g)\to R_m(\gsh),$$
dual to the 'edge homomorphism' $J_m(\gsh)\to J_m(\g)$.\par In the
case of a Lie algebra $\g$, a local homomorphism $$R_m(\g)\to
(S,\m)$$ to a local artin algebra is given by a 'morphic' (i.e.
comultiplicative) cochain in $H^0(J_m(\g\otimes\m))$, which has the
form
$$\epsilon(\phi)=(\phi, ...,\frac{1}{m!}\phi^{\times m})$$ where
$\phi\in\g^1\otimes\m$ satisfies the Kodaira-Spencer integrability
condition\nl
$$\bp\phi=\frac{-1}{2}[\phi,\phi].\tag ??$$ A similar fact holds for a Lie
atom $\gsh$: a local homomorphism $R_m(\gsh)\to (S,\m)$ is given by
a morphic cochain $$\epsilon (\phi,\psi)=(\epsilon(\phi), (\psi,
\epsilon(\phi)\times\psi), (\psi^2,
\epsilon(\phi)\times\psi^2),...)\in H^0(J_m(\gsh\otimes\m))$$
satisfying the integrability conditions
$$\bp\phi=\frac{-1}{2}[\phi,\phi],$$
$$\bp\psi+i(\phi)=\frac{-1}{2}\langle\phi,\psi\rangle.\tag ??$$
Similarly when $\gsh$ is replaced by $\g^+.$ Now in between
$H^0(J_m(\gsh\otimes\m))$ and\nl $H^0(J_m(\g^+\otimes\m))$ is
$H^0(J_m(\g^+_\m))$. Note that a morphic cochain in
$H^0(J_m(\g^+_\m))$, i.e. a cochain $\epsilon=\epsilon(\phi,\psi)$
as above with $\phi\in\g^0\otimes\m, \psi\in\h^+_\m$ and satisfying
integrability, induces  for $i=1,...,e$ a cochain in $H^0(J_m(\g\sbr
i._\m)), $ hence another morphic cochain $\epsilon_i\in
H^0(J_m(\g\sbr i.\otimes(\m/\m^{i+1})))$. Clearly the image of
$\epsilon$ in $H^0(J_m(\g^+\otimes\m))$ 'is' (i.e. comes from)
$\epsilon_e$. Thus we get a compatible set of homomorphisms
$$\matrix R_m(\g^+)&\to&R_m(\g\sbr e.)&\to...\to&
R_m(\g\sbr i.)&\overset{i_i^*}\to{\to}...\to &R_m(\gsh)=R_m(\g\sbr 1.)\\
&\alpha\searrow&\alpha_e\downarrow&&\alpha_i\downarrow &&\alpha_1
\downarrow
 \\
&&S&\to...\to&S/\m^{i+1}&\to...\to &S/\m^2\endmatrix\tag ??$$ and
$i_i^*$ is surjective by ??. Of course if
 $\gsh$ is a Lie pair then
 $$R_m(\g^+)=...=R_m(\g^i)=...=R_m(\gsh).$$ in general,
 clearly $\alpha_i$ factors through $R_i(\g\sbr i.).$
 Moreover, it follows easily from ?? that given an element,
 comultiplicative or not, $\epsilon\in H^0(\g^+\otimes\m)$,
 $\epsilon$ is in the subspace $H^0(\g^+_\m)$ iff the induced map
 $R_m(\g^+)\to S$ admits a factorization as in ??.
 There is no loss of generality in assuming $\alpha_e$ is
 surjective. In this case, all diagrams (??) are obtained as
 follows. Write
 $$R_m^+=R_m(\g^+),K=\ker(R_m^+\to S), J_i=\ker(R_m^+\to
 R_m(\g\sbr i.)).$$ Then the conditions on $K$ are
 $$\m_R^{e+1}\subseteq K, J_i\subseteq K+\m_R^{i+1}, i=1,...,e.$$
 Geometrically, this means $S$ corresponds to a germ of a variety
 that is tangent to order $i$ to the locus $\Spec(R_m(\g\sbr
 i.))\subset\Spec (R_m^+).$\par
 Note that if we have a global family over a smooth base $B$ that
 yields a $\g^+$-deformation at a general point of $B$, then $B$ is
 endowed at its general point with a chain of distributions
 corresponding to the subatoms $\g^i.$ If the deformation is in fact
 a $\gsh$-deformation, then the distribution corresponding to
 $\g\sbr 1.=\gsh$ already coincides with the full tangent space to
 $B$, and it follows that in ??, the map $\alpha_e$ already factors
 through a map $R_m(\gsh)\to S.$

  \remark{Remark}If $\h$ happens to be a Lie algebra and
$i$ a Lie homomorphism we may think of $R_m(\g/\h)$ as a formal
functorial substitute for the fibre of the induced homomorphism
$R_m(\h)\to R_m(\g)$, i.e. $R_m(\g)/\m_{R_m(\h)}. R_m(\g)$. It is
important to note that this fibre involves only the $\g-$module
structure on $\h$ and not the full algebra structure.
\endremark
Now we may associate to a $\gsh-$deformation a
$J_m(\gsh_\m)-$cocycle as follows. In terms of modules, we have, as
in ??, for each let $\gsh$ module $E_1<E_2$ and right module
$F_1>F_2$ a dual pair of $\h-$ maps
$$A: E_2\otimes S\to E_2^\phi,$$
$$B: F_1^\phi\to F_1\otimes S.$$
As usual we represent $E_i^\phi$ by a resolution of the from
$(E_i^.\otimes S, \bp +\phi), i=2,3$ and $E_i\otimes S$ by
$(E_i^.\otimes S, \bp);$ similarly for $F.$. Then the maps $A,  B$
can be represented simultaneously in the form $\exp(\pm v_\alpha),
v_\alpha\in \h^{+0}\m,$ and we get a pair of commutative diagrams
$$\matrix E_2^i\otimes S&\overset\bp\to\to&E_2^{i+1}\otimes S\\
\exp(v_\alpha)\downarrow&&\downarrow \exp(v_\alpha)\\
E_2\otimes S&\overset{\bp+\phi}\to\to&E_2^{i+1}\otimes S
\endmatrix\tag 3.3$$
$$\matrix
F_1^i\otimes S&\overset{\bp+\phi}\to\to&F_1^{i+1}\otimes S\\
\exp(-v_\alpha)\downarrow&&\downarrow \exp(-v_\alpha)\\
F_1\otimes S&\overset{\bp}\to\to&F_1^{i+1}\otimes S\endmatrix
$$
Universally, in terms of e\v{C}ech torsor data
$(\Psi_{\alpha\beta}=\exp (u_{\alpha\beta}), v_\alpha)$ for a good
atlas as above, the condition that the $\exp(v_\alpha)$ should glue
together to a globally defined map left $\h-$ map $A$ is
$$\Psi_{\alpha\beta}\circ \exp(v_\alpha )=\exp(v_\beta).\tag 2.3$$
By analogy with the procedure of \cite{R, p.61}, this condition may
be analyzed in terms of Kodaira-Spencer cochains, as follows. Write
$$\Psi_{ab}=\exp(s_\alpha)\exp(-s_\beta),$$
$$s_\alpha\in\g^0(U_\alpha)\otimes\m$$
where $\g^.$ is a suitable acyclic (soft or injective) \dgla
resolution of $\g.$ Thus
$$\exp(-s_\alpha)\exp(v_\alpha)=\exp(-s_\beta)\exp(v_\beta)$$ so
this element is globally defined and can be written as
$$\exp(\psi),\ \  \psi\in\Gamma(\h^0)^+_\m,$$
 where $\h^.$ is a suitable
$\g^.$-module resolving $\h$ and $\h^{.+}$ is a suitable \dgla
resolution of $\h^+$. As observed in \cite{Rcid},
$$\phi=\exp(-s_\alpha)d\exp(s_\alpha)\in\g^1\otimes\m$$ is globally
defined independent of $\alpha$ and is the Kodaira-Spencer cochain
defining the deformation. Then a formal calculation in
$\U_S(\Gamma(\h^{.})^+_\m)$ shows that
$$d\exp(\psi)=-\phi\exp(\psi),d\exp(-\psi)=\exp(-\psi)\phi$$ which
formally implies $$d\psi+i(\phi)=\frac{-1}{2}\langle
\phi,\psi\rangle $$

Now we may associate to a $\gsh-$deformation a $J_m(\gsh)-$cocycle
as follows. By definition, we have a pair of pairs of
$\g-$deformations with a dual pair of $\h-$ maps
$$A: E_1\otimes S\to E_2^\phi,$$
$$B: E_3^\phi\to E_4\otimes S.$$
As usual we represent $E_i^\phi$ by a resolution of the from
$(E_i^.\otimes S, \bp +\phi), i=2,3$ and $E_i\otimes S$ by
$(E_i^.\otimes S, \bp).$ Then the maps $A,  B$ can be represented
simultaneously in the form $j\pm v, v\in \h^0\otimes\m_S,$ and we
get a pair of commutative diagrams
$$\matrix E_1^i\otimes S&\overset\bp\to\to&E_1^{i+1}\otimes S\\
j+v\downarrow&&\downarrow j+v\\
E_2\otimes S&\overset{\bp+\phi}\to\to&E_2^{i+1}\otimes S
\endmatrix\tag 3.3$$
$$\matrix
E_3^i\otimes S&\overset{\bp+\phi}\to\to&E_3^{i+1}\otimes S\\
j-v\downarrow&&\downarrow j-v\\
E_4\otimes S&\overset{\bp}\to\to&E_4^{i+1}\otimes S
\endmatrix\tag 3.4$$
whose commutativity amounts to a pair of identities in $\U(\g,\h)$:
$$\phi .v =-\bp (v) -i(\phi),\ \ v.\phi=\bp(v)+i(\phi)\tag 3.5$$
which imply
$$\bp(v)+i(\phi)=-\1/2 \langle  \phi,v\rangle .\tag 3.6$$
Now the identity (3.6) together with the usual integrability
condition
$$\bp(\phi)=-\1/2 [\phi,\phi]\tag 3.7$$
imply that, setting
$$\epsilon (\phi)=(\phi, \phi\times\phi, ...),$$
$$\epsilon(\phi, v)=
(v, \phi\times v, \phi\times\phi\times v, ...),$$ the pair
$$\eta (\phi, v):=( \epsilon (\phi), \epsilon (\phi, v))$$
constitute a cocycle for the complex $J_m(\gsh)\otimes\m_S.$ This
cocycle is obviously 'morphic' or comultiplicative, hence gives rise
to a ring homomorphism
$$\alpha^\sharp=\alpha^\sharp(\phi,v) : R_m(\gsh)\to S\tag 3.8 $$
lifting the usual Kodaira-Spencer homomorphism $\alpha(\phi):
R_m(\g)\to S$.\par Conversely, given a homomorphism $\alpha^\sharp$
as above, with $S$ an arbitrary artin local algebra, clearly we may
represent $\alpha^\sharp$ in the form $\eta (\phi, v)$ as above
where $\phi$ and $v$ satisfy the conditions (3.6) and (3.7). Then in
the enveloping algebra $\U(\h^+)$ we get the identity
$$i(\phi)=-\bp(v)-\1/2 [i(\phi), v].\tag 3.9$$
Plugging the identity (3.9) back into itself we get, recursively,
$$i(\phi)=-\sum\limits_{k=0}^\infty \frac{(-1)^k}{2^k}
\ad (v)^k(\bp(v))\tag 3.10$$ (the sum is finite because $\m_S$ is
nilpotent), from which the identities 3.5 follow formally. Hence the
diagrams 3.3 and 3.4 commute, so we get a $ \gsh-$deformation
lifting $\phi.$\par In particular, applying this construction to the
identity element of $S=R_m(\gsh)$, thought of as an element of
$\H^0(J_m(\gsh))\otimes\m_S$, we obtain a 'tautological'
$\gsh-$deformation which we denote by
$$u_m^\sharp=(\phi_m,v_m)$$
and we get the following analogue (and extension) of Theorem 3.1:
\proclaim{Theorem 3.2} Let $\gsh$ be an admissible \dg Lie atom such
that $\g$ has central sections. Then for any $\gsh-$ deformation
$(\phi,v)$ over an artin local algebra $S$ of exponent $m$ there
exists a nonempty  ${\text{Aut}}(\phi,v)$-torsor of isomorphisms of
$\gsh$-deformations
$$(\phi,v)\simeq \alpha^\sharp(\phi,v)^*(u_m^\sharp).\qed$$
\endproclaim

\endcomment
\remark{Examples 1.1.4, conclusion}
Applying Theorem 3.3 to our standard examples, we see the following.
\item{\bf{A-B}} When $E_1<E_2$
are coherent sheaves on a  complex projective scheme (or vector
bundles on a compact complex manifold) $X$, the universal
$\gl(E_1<E_2)$-deformation just constructed is the formal completion
of the quot scheme Quot$(E_2)$ at the point $(E_2\to E_2/E_1).$ The
projectivity or compactness hypothesis $X$  is only needed to ensure
admissibility. See a forthcoming article \cite{Rsela} for applications
of Lie atoms to local moduli and parameter spaces for schemes.
\item{\bf C} For heat deformations, we have constructed the
universal heat deformation of $E$ (again over a compact complex
manifold, to ensure admissibility).\item {\bf {D}} The universal
$N_{Y/X}[-1]$-deformation is the formal germ of the Hilbert scheme
or Douady space, first constructed by Kodaira \cite{K}. Note
that because $N_{Y/X}[-1]$ is (cohomologically) supported on $Y$,
only compactness of $Y$ is needed to ensure admissibility. For a
holomorphic map $f:Y\to X$, the universal $T_f$-deformation is the
universal deformation of the map $f$, first constructed by other
means by Horikawa \cite{Ho}. When $f$ is generically immersive, the
universal $N_{f/X}[-1]$-deformation is the universal deformation of
$f$ with fixed $X$, and similarly for $N_{f/X}$ when $f$ is
generically submersive.\item{\bf E} The universal $T_X\otimes \O_Y[-1]$- deformation is the universal deformation of the map $Y\to X$ fixing $X$ and $Y$.\par
One consequence of applying Theorem 3.3 on these examples is that Corollary 1.2.3 applies to them, and we conclude \proclaim{Corollary 3.4} In each of the above examples {\bf {A-E}}, if the relevant Lie atom $\gsh$ satisfies $$h^1(\gsh)-h^2(\gsh)>0,$$ then the corresponding object admits a nontrivial deformation.\endproclaim As one example, if a submanifold $Y\subset X$ satisfies $(h^0-h^1)(N_{Y/X})>0$, then $Y$ moves in $X$.

 \comment\item{\bf F} We
conclude with a new example generalizing the previous 2. Let $f:Y\to
X$ be a holomorphic map of complex manifolds. In \cite{Rdm} we
introduced a Grothendieck topology $T(f)$, whose open sets are pairs
of opens $(U\subset Y, V\subset X)$ such that $f(U)\subset V.$  On
$T(f)$ we can define a sheaf $T_f$ of Lie algebras whose sections on
$(U,V)$ consist of the pairs of vector fields $u\in T_Y(U), v\in
T_X(V)$ such that $f^*(v)=df(u)$ where $df:TY\to f^*T_X$ is the
differential. There is an obvious map $T_f\to f^*T_X$ which is
injective if $f$ is generically injective, and an obvious
$T_f-$module structure on $f^*T_X$, making the pair $(T_f, T_X)$ a
sheaf of Lie atoms on $T(f)$, which is admissible if $X,Y$ are
compact. Thus Theorem 3.2  gives a Lie-atomic construction of the
universal deformation of the map $f$, first constructed by other
means by Horikawa \cite{Ho}.

\endcomment
\endremark
\remark{Remark}{\rm After this was written, the author became aware
of a preprint '$L_\infty$ structures on mapping cones'
(arxiv:math/QA/0601312) by D. Fiorenza and M. Manetti which, quoting
an earlier version of this paper, considers related problems (in the
case of Lie pairs only). Fiorenza and Manetti consider
'deformations' in a purely formal, Lie-theoretic sense, as solutions
of the Maurer-Cartan integrability equation, ignoring the connection
with the geometric view of deformations as torsors and questions of
existence of universal deformations, and therefore also the subtle
question of the relation between formal or cohomological equivalence
of Maurer-Cartan solutions (measured e.g. via the Jacobi-Bernoulli
complex or something similar), and equivalence or conjugacy of
deformations, either viewed as torsors or realized on a particular
module (compare the Introduction to \S2 above). Also, their results
do not apply to general Lie atoms, hence they cannot be applied to
Heat deformations. The purely formal results of Fiorenza-Manetti
follow easily from our technique of Jacobi-Bernoulli complex. On the
other hand the our main results here do focus on geometric
deformations, especially universal ones, and for that reason
involve, of necessity, some restrictive hypothesis such as
finiteness and asymmetry.\par As mentioned above, the methods of
this paper, such as Jacobi-Bernoulli cohomology, have now been
extended from Lie atoms to Semi-simplicial Lie algebras (SELA). It
can also be shown that one can associate a SELA $\Cal T_X$ to any
algebraic scheme $X/\C$, so that deformations of $X$ as scheme can
be expressed in terms of the Jacobi-Bernoulli cohomology of $\Cal
T_X$. See \cite{Rsela} for details.\par
I am grateful to Sergei Merkulov for pointing out an inaccuracy
in the formulation of (1.2.4); see his work now in progress
for applications of JB ideas in the context of operads.
\endremark
\heading References\endheading
\item{[At]}
Atiyah, M.F.: Complex analytic connections in fibre bundles.
Transac. Amer. Math. Soc. 85 (1957), 181-207.
\item{[Hi]} Hitchin, N. J.: Flat connections and geometric
quantisation. Commun. math. Phys. 131, 347-380 (1990).
\item{[Ho]} Horikawa, E.: Deformations of holomorphic maps, I
(J. Math. Soc. Japan 25 (1973), 372-396), II (ibid. 26 (1974),
647-667), III (Math. Ann. 222 (1976), 276-282).
\item{[K]} Kodaira, K.: A theorem of completeness of characteristic
systems for analytic families of compact submanifolds of complex
mainifolds. Ann. Math. 75 (1962), 146-162.
\item{[Kol]} Koll\'ar, J.: Rational curves on algebraic varieties. Springer.
\item{[Rcid]}
Ran, Z.: Canonical infinitesimal deformations. J. Alg. Geometry 9,
43-69 (2000).
\item{[Ruvhs]}Ran, Z.: Universal variations of Hodge structure.
Invent. math. 138, 425-449 (1999).
\item{[Rrel2]} Ran, Z. : Jacobi cohomology, local geometry of moduli
spaces, and Hitchin connections. Proc. Lond. Math
Soc. (3) 92 (2006), 545-580.
\item{[Rsela]} Ran, Z. : Deformations of schemes and maps (in
preparation).

\item{[VLP]} Verdier, J.-L., Le Potier, J., edd.: Module
des fibr\'es stables sur les courbes alg\'ebriques. Progr. in Math.
54 (1985), Birkh\"auser.
\item{[W]}
Welters, G.: Polarized abelian varieties and the heat equations.
Compos. math. 49, 173-194 (1983).

\enddocument